\documentclass{amsart}

\usepackage{amsmath,amsthm,amssymb}
\usepackage{amsfonts}
\usepackage{rotating}
\usepackage{euscript}
\usepackage{pst-node}
\usepackage{epsfig}

\def\be{\begin{equation}}
\def\ee{\end{equation}}

\def\C{{\mathbb C}} 
 
\def\P{{\mathbb P}}
\def\Z{{\mathbb Z}}
\def\R{{\mathbb R}} 
\def\Q{{\mathbb Q}}
\def\ord{{\rm ord}}

\def\e{\eqref}
\def\phi{{\varphi}}
\def\v{{\varepsilon}} 
\def\deg{{\rm deg\,}}

\def\Log{{\rm Log\,}} 
\def\Gal{{\rm Gal\,}}

\def\cos{{\rm cos\,}}

\def\mod{{\rm mod\ }}
\def\qed{$\ \ \Box$ \vskip 0.2cm}

\def\bp{\begin{proposition}}
\def\ep{\end{proposition}}

\def\bt{\begin{theorem}}
\def\et{\end{theorem}}
\def\br{\begin{remark}}
\def\er{\end{remark}}
\def\be{\begin{equation}}
\def\bee{\begin{equation*}}
\def\la{\label}
\def\l{\label}

\def\tilde{\widetilde}
\def\hat{\widehat}
\def\ee{\end{equation}}
\def\eee{\end{equation*}}
\def\bl{\begin{lemma}}
\def\el{\end{lemma}}
\def\br{\begin{remark}}
\def\er{\end{remark}}
\def\bc{\begin{corollary}}
\def\ec{\end{corollary}}
\def\pr{\noindent{\it Proof. }}

\def\bd{\begin{definition}}
\def\ed{\end{definition}}
\input epsf.sty

\newtheorem{theorem}{Theorem}[section]
\newtheorem{lemma}[theorem]{Lemma}
\newtheorem{corollary}[theorem]{Corollary}
\newtheorem{proposition}[theorem]{Proposition}
\newtheorem{remark}{Remark}[section]

\begin{document}
\title[On functions orthogonal to all powers of a given function]{On rational functions orthogonal to all powers of a given rational function 
on a curve}
\author{F. Pakovich}
\address{Department of Mathematics, Ben-Gurion University
of the Negev, P.O.B. 653, Beer-Sheva, Israel} 
\email{pakovich@cs.bgu.ac.il}
\date{}

\begin{abstract} In this paper we study the generating function $f(t)$ for the sequence of the moments $\int_{\gamma}P^i(z)q(z)d z,$ $i \geq 0,$ where $P(z),q(z)$ are rational functions of one complex variable and 
$\gamma$ is a curve in $\C.$ We calculate an analytical expression for $f(t)$ and provide conditions implying the rationality and the vanishing of $f(t).$  
In particular, for $P(z)$ in generic position we give an explicit 
criterion for a function $q(z)$ to be orthogonal to all powers of $P(z)$ on  
$\gamma.$ Besides, we 
prove a stronger form of the Wermer theorem, describing 
analytic functions satisfying $\int_{S^1}h^i(z)g^j(z)g'(z)d z=0,$ $i\geq 0,$ $j\geq 0,$ in the case where the functions $h(z),g(z)$ are rational.
We also generalize the theorem of Duistermaat and van der Kallen about Laurent polynomials $L(z)$ whose integral positive powers have no constant term, and prove other results about  
Laurent polynomials $L(z),m(z)$ satisfying $\int_{S^1}L^i(z)m(z)d z=0,$ $i\geq i_0.$

\end{abstract}

\maketitle

\section{Introduction}
In this paper we study the generating function $f(t)=\sum_{i=0}^{\infty}m_it^i$ for the sequence of the moments
\be \la{1} m_i=\int_{\gamma}P^i(z)q(z)d z, \ \ \ i\geq 0, \ee
where $P(z)$, $q(z)$
are rational functions of one complex variable and $\gamma$ is a curve in $\C$. In particular, we study conditions under which $f(t)$ is a rational function, a polynomial, or an identical zero.

The main motivation for such a study is the relation of the function $f(t)$ with the classical Center-Focus problem of Poincar\'e. Let $F(x,y),$ $G(x,y)$ be real-valued functions of $x,y$ analytical in a neighborhood of the origin in $\R^2$ and vanishing at the origin together with their first derivatives. The Poincar\'e problem is to find conditions
under which all solutions of the system 
\be \la{sys1}
\left\{
\begin{array}{rcl}
\dot{x} & =& -y+F(x,y), \\ 
\dot{y} & =&x+G(x,y),\\
\end{array} 
\right.
\ee
around zero are closed (see e.g. the recent survey \cite{c1} and the bibliography therein).
Despite the efforts of many researchers 
this problem remains open even in the case where 
$F(x,y)$ and $G(x,y)$ are polynomials of degree 3, 
and any advances in its understanding are of a great interest.

It was shown in \cite{cher} that if $F(x,y),$ $G(x,y)$ are homogeneous polynomials of the same degree, then 
one can construct another polynomials $f(x,y),$ $g(x,y)$ such that \eqref{sys1} has a center if and only if all
solutions of the trigonometric Abel equation 
\be \la{sys2}
\frac{d r}{d \phi}=f(\cos \phi,\sin \phi)\, r^2+g(\cos \phi,\sin \phi)\, 
r^3 \ee
with $r(0)$ small enough
are periodic on $[0,2 \pi]$. 
In its turn, the trigonometric Abel equation can be transformed by an exponential
substitution into the equation
\be \la{sys3}
\frac{d y}{d z}=l(z)y^2 + m(z)y^3, 
\ee
where $l(z)$ and $m(z)$ are Laurent polynomials. 
Furthermore, all solutions of \eqref{sys2} 
with $r(0)$ small enough are periodic on $[0,2 \pi]$ if and only if all solutions of \eqref{sys3} with $y(1)$ small enough  
are non-ramified along $S^1$. 
By the analogy with system \eqref{sys1}, if 
all solutions of \eqref{sys2} 
with $r(0)$ small enough are periodic on $[0,2 \pi]$ (resp. if 
all solutions of \eqref{sys3} with $y(1)$ small enough  
are non-ramified along $S^1$) we will say that
Abel equation \eqref{sys2} (resp. Abel equation \eqref{sys3}) has a center. 

In the series of papers \cite{bfy1}-\cite{bfy4} the following modification of the center problem for equation \eqref{sys3} was proposed:
find conditions under which 
for any solution $y(z)$ of the Abel differential equation
\be \la{sys4}
\frac{d y}{d z}=p(z)y^2+q(z)y^3
\ee with {\it polynomial} coefficients $p(z),$ $q(z)$
the equality $y(1)=y(0)$ holds whenever $y(0)$ is small enough (as above, in case if this condition is satisfied we will say 
that Abel equation \eqref{sys4} has a center). 
This modification seems to be easier than the initial problem at the same time  
keeping its main features, and in this context can be considered as a
simplified form of the classical Center-Focus problem of Poincar\'e.

The center problem 
for Abel equation \eqref{sys4}
naturally leads to the following 
``polynomial moment problem'': find conditions under which 
polynomials $P(z),$ $q(z)$ satisfy the system  
\be \la{11} \int_{0}^1P^i(z)q(z)d z=0, \ \ \ i\geq 0, \ee
or, in other words, find conditions implying the vanishing of the corresponding function $f(t)$.
The center problem for Abel equation \eqref{sys4}
is related to the 
polynomial moment problem in several ways.
For example, it was shown in \cite{bfy3}
that for the parametric version  
$$
\frac{d y}{d z}=p(z)y^2+\varepsilon q(z)y^3
$$
of \eqref{sys4} the ``infinitesimal'' center condition with respect to $\varepsilon$ 
reduces to \eqref{11} with $P(z)=\int p(z) d z.$
On the other hand, it was shown in \cite{bry} that ``at infinity'' (under an appropriate projectivization of the 
parameter space) the center condition for \eqref{sys4} also reduces to \eqref{11}.
For other results relating 
the center problem for equation \eqref{sys4} and the polynomial moment problem see \cite{bry} and
the bibliography therein.

Posed initially as an intermediate step on the way to the solution of the Poincar\'e
problem, the polynomial moment problem turned out to be quite subtle question unexpectedly related with such 
branches of mathematics as the Galois theory and representations of groups. This problem has been studied in the papers \cite{bfy1}-\cite{bfy5}, \cite{c}, \cite{pa1}-\cite{pm}, \cite{pakk}, \cite{ro}
with the final solution achieved in \cite{pm}, \cite{pakk}. 
Before formulating the results of \cite{pm}, \cite{pakk} explicitly let us introduce the following  
``composition condition'' which being imposed on $P(z)$ and $Q(z)=\int q(z)d z$
implies both equalities \eqref{11} and a center for \eqref{sys3}.
Namely, suppose that
there exist polynomials $\tilde P(z),$ $\tilde Q(z),$ $W(z)$ such that
the  equalities \be \l{2}
P(z)=\tilde P(W(z)), \ \ \ 
\ Q(z)=\tilde Q(W(z)), \ \ \ W(0)=W(1)
\ee hold.
Then changing the variable $z$ to $W(z)$
and taking into account that for any path $\gamma$ connecting $0$ and $1$ the path $W(\gamma)$ is closed, it is easy to see that equalities \eqref{11} hold and Abel equation 
\eqref{sys3} has a center.
Furthermore, the main conjecture concerning the center problem for the Abel equation
(``the composition conjecture for the Abel equation") states 
that \eqref{sys3} has a center if and only \eqref{2} holds.

For many classes of $P(z)$ the composition condition \eqref{2} turns 
out to be necessary for equalities \eqref{11} to be satisfied.  
For instance, this is true  
if $0,1$ are not critical points of $P(z)$   
(\cite{c}) or if $P(z)$ is indecomposable that is can not be represented as a composition of two non-linear polynomials (\cite{pa2}).
Nevertheless, 
this is not true in general \cite{pa1}, and a right description of 
solutions of \eqref{11} is as 
follows \cite{pm}:
non-zero polynomials $P(z),$ $q(z)$ satisfy \eqref{11} 
if and only if $Q(z)=\int q(z)d z$ can be represented as a {\it sum} of polynomials $Q_j(z)$ such that
\be \l{cc}
P(z)=\tilde P_j(W_j(z)), \ \ \
Q_j(z)=\tilde Q_j(W_j(z)), \ \ \ {\rm and} \ \ \ W_j(0)=W_j(1)
\ee
for some polynomials $\tilde P_j(z), \tilde Q_j(z), W_j(z)$.
Moreover, it was shown in \cite{pakk} that actually any solution of the polynomial moment problem may be obtained as a sum of at most two reducible solutions and that the corresponding reducible solutions may be described in a
very explicit form.

In the same way as the center problem for Abel equation \eqref{sys4}
leads to the polynomial moment problem, the center problem 
for equation \eqref{sys3} leads to the following 
``Laurent polynomial moment problem", which is the main motivation for investigations of this paper:
describe Laurent polynomials $L(z), m(z)$ such that
\be\la{lloo}
\int_{S^1} L^i(z)m(z) d z= 0,\ \ \  i\geq 0.
\ee 
An analogue of condition \eqref{2} in this setting
is that
there exist a Laurent 
polynomial $W(z)$ and polynomials $\tilde L(z),$ $\tilde m(z)$ such that the equalities 
\be \la{uxx} L(z)=\tilde L(W(z)), \ \ \ m(z)=\tilde m(W(z))W'(z)\ee
hold. Clearly, \eqref{uxx} implies \eqref{lloo}. Furthermore, if for given $L(z)$ 
there exist several such $m(z)$, then \eqref{lloo} is satisfied for their sum. 
However, in distinction with the polynomial moment problem other mechanisms for \eqref{lloo} to be satisfied also exist. For example, if 
$L(z)=\tilde L(z^d)$ for some $d>1$, then the residue calculation shows that condition
\eqref{lloo} is satisfied for any Laurent polynomial containing no terms $cz^{n}$ with $n\equiv -1\, \mod d.$ 

In addition to topics related to the Poincar\'e problem,
the questions concerning 
the function $f(t)$ appear also in other circumstances. Let us mention for example the following particular case  
of the Mathieu conjecture concerning compact Lie groups \cite{mat}
proved by Duistermaat and van der Kallen \cite{dui}:
if all integral positive powers of a Laurent polynomial $L(z)$ have no constant term, then
$L(z)$ is either a polynomial in $z$, or a polynomial in $1/z.$
Clearly, the assumption of this theorem is equivalent to the condition that the function
$f(t)$ for $P(z)=L(z)$, $m(z)=1/z,$ and $\gamma=S^1$ 
is a constant. Observe a delicate difference with the Laurent polynomial moment problem: in the assumption 
of the theorem of Duistermaat and van der Kallen
it is not required that integral in \eqref{lloo} equals zero for $i=0$. Furthermore, the addition of such a condition would make the problem meaningless since for $i=0$ and $m(z)=1/z$ integral in \eqref{lloo} is distinct from zero for any non-zero $L(z).$ Notice also that the question about description of polynomials $P(z),$ $q(z)$ for which the function $f(t)$ is equal to a constant, and not just to zero as in the polynomial 
moment problem, was recently raised by Zhao \cite{gene} in relation with his conjecture about images of commuting differential operators 
(such a description is not an immediate 
corollary of the solution of the polynomial moment problem given in  \cite{pm}). 

Finally, observe that the classical Wermer theorem \cite{w1}, \cite{w2}, describing 
analytic functions on $S^1$ 
satisfying \be \la{weeer} \int_{S^1}h^i(z)g^j(z)g'(z)d z=0, \ \ \ i,j\geq 0,\ee in the case where the functions $h(z),g(z)$ are rational, obviously also is related  
to the subject of this paper. Here, however, the assumption 
requires that an {\it infinite} number of functions $f_j(t),$ $j\geq 0,$ corresponding to $P_j(z)=h(z),$ $q_j(z)=g^j(z)g'(z),$ and $\gamma=S^1,$ vanish simultaneously. Notice that for rational $h(z)$ and $g(z)$ the Wermer theorem is equivalent to the following statement \cite{bry0}: 
a necessary and sufficient condition for $h(z),$ $g(z)$ to satisfy 
equalities \eqref{weeer} is that there exist rational functions 
$\tilde h(z),$ $\tilde g(z),$ $w(z)$ such that
\be \la{werr} 
h(z)=\tilde h(w(z)), \ \ \ 
\ g(z)=\tilde g(w(z)),
\ee and the curve $w(S^1)$ is homologous to zero in $\C\P^1$ with 
poles of $\tilde f(z)$ and $\tilde g(z)$ removed.

This paper is organized as follows. In the second section we calculate an analytical expression for $f(t).$ Our approach here is similar to the one of \cite{pry}; however, in contrast to \cite{pry}, we
obtain an explicit  
analytical expression for $f(t)$.  
Notice that formulas obtained imply in particular that if $\gamma$ is closed, then $f(t)$ is an algebraic function from the field $K_P$ generated over $\C(z)$ 
by the branches $P^{-1}_{i}(z),$ $1\leq i \leq n,$ of the algebraic function $P^{-1}(z)$ inverse to $P(z)$, 
while if $\gamma$ is non-closed, then $f(t)$ is a linear combination 
of branches of the logarithm with coefficients from $K_P$. 
Another corollary is that if $\gamma$ is closed and 
homologous to zero in $\C\P^1$ with 
poles of $P(z)$ removed, then the function $f(t)$ is rational
for {\it any} $q(z)$.
Following \cite{pry}, if $P(z)$ and $\gamma$ satisfy 
the last condition we will say that {\it poles of $P(z)$ lie on one side of $\gamma$}.

Although the analytical expression for $f(t)$ obtained in the second section is explicit, in general this expression by itself does not allow us to conclude whether $f(z)$ is 
rational or vanishes identically, and the third section of the paper 
is devoted to these questions.  
We show that if poles of $P(z)$ do not lie on one side of $\gamma$, then 
$f(t)$ is rational if and only if the superpositions of the rational function 
$(q/P')(z)$ with
branches of $P^{-1}(z)$ satisfy 
a certain system of equations 
\be \la{sisa}
\sum_{i=1}^nf_{s,i}\left(\frac{q}{P'}\right)(P^{-1}_{i}(z))=0, \ \ \ \ \ \ f_{s,i} \in \Z,  \ \ \ \ \ \ 1\leq s\leq k, 
\ee where $f_{s,i}$ and $k$ depend on $P(z)$ and $\gamma$ only. 
This result generalizes the corresponding criterion for polynomials given in \cite{pp} and relies on the same ideas. In particular, we use combinatorial objects called ``constellations" (similar to what is called ``Dessins d'enfants") which represent the monodromy group $G_P$ of the algebraic function  $P^{-1}(z)$ in a combinatorial way.

In the third section we also show that if $P(z)$ and $q(z)$ satisfy the conditions  
\be \la{condi} q^{-1}\{\infty\}\subseteq P^{-1}\{\infty\},\ \ \ P(\infty)=\infty,\ee then the rationality of $f(t)$ yields that $f(t)\equiv 0.$ Notice that this result implies immediately the theorem of Duistermaat and van der Kallen 
cited above. Indeed, if $L(z)$ is not a polynomial in $z$ or in $1/z$, then for $P(z)=L(z),$ $q(z)=1/z$ conditions \eqref{condi} 
are satisfied and therefore the equality
$f(t)=c$ would imply that $c=0$ in contradiction with the fact that 
for $i=0$ the integral in \eqref{lloo} is not zero.  Besides, the mentioned result implies the following statement which gives the answer to the question of Zhao: if integrals in \eqref{11} vanish for all $i\geq i_0$, 
then they vanish for all $i\geq 0$ and therefore the polynomial $Q(z)$ is a sum of polynomials $Q_j(z)$ such that \eqref{cc} holds.

In the fourth section, using system \eqref{sisa} 
and the characteristic property of permutational matrix representations of doubly transitive groups, 
we show that if $P(z)$ is in generic
position and poles of $P(z)$ do not lie on one side of $\gamma$, then the function $f(t)$ is rational if and only if 
the function $q(z)$ has the form  \be \la{fif} q(z)=\tilde q(P(z))P'(z)\ee for some rational function $\tilde q(z).$
Besides, we show that for $P(z)$ as above a rational function $q(z)$ 
satisfies 
$$\int_{\gamma}P^i(z)q(z)d z=0, \ \ \ i\geq 0,$$ if and only if \eqref{fif} holds for some rational function $\tilde q(z)$ whose
poles lie outside the curve $P(\gamma).$

In the fifth section we study conditions for 
vanishing of double moments \eqref{weeer} for rational $h(z),$ $g(z)$ in a more general setting where the integration path is an arbitrary
curve 
and the conditions $i\geq 0,$ $j\geq 0$ may be weaker. We prove two results which precise and generalize the Wermer theorem in the special case where $h(z),$ $g(z)$ are rational. 

In the sixth section we obtain some preliminary results concerning the Laurent polynomial moment problem. In particular, we deduce from the results of the fourth section its solution in the case where $L(z)$ is in generic position. Besides, we prove two generalizations of the theorem of Duistermaat and van der Kallen.

In the seventh section we study the following problem: how many 
first integrals in \eqref{lloo} should vanish in order to conclude that 
all of them vanish. Using the fact that the corresponding function 
$f(t)$ is contained in the field $K_L$, we give a bound which depends on degrees of Laurent polynomial $L(z)$ and $m(z)$ only.

Finally, the eighth section is devoted to relations between the condition 
that $f(t)$ is rational and the condition that there exist rational functions $\tilde q(z),$ $\tilde P(z)$,
$W(z)$ such that $\deg W(z)>1$ and 
\be \la{qaz} P(z)=\tilde P(W(z)), \ \ \ q(z)=\tilde q(W(z))W'(z)\ee
(notice that \eqref{fif} is a particular case of \eqref{qaz} where $W(z)=P(z)$ and $\tilde P(z)=z$, while 
\eqref{2}, without the requirement $W(0)=W(1)$, is equivalent to \eqref{qaz} for $q(z)=Q'(z)$).
More precisely, we are interested in conditions which imply that for 
given $P(z)$ and $\gamma$ such that poles of $P(z)$ do not lie on one side of $\gamma$, the rationality of $f(t)$ implies that \eqref{qaz} holds
for some  $\tilde q(z),$ $\tilde P(z)$, and $W(z)$ with $\deg W(z)>1$. Using 
a general algebraic result of Girstmair \cite{gi2} we give such a criterion and discuss an
explanatory example.

\section{Analytic expression for $I_{\infty}(t)$}
\subsection{Definition of the function $I_{\infty}(t)$}
Let $P(z),q(z)$ be rational functions and $\gamma\subset \C$ be a curve. 
We always will assume that $\deg P(z)>0$ and $q(z)\not\equiv 0.$ We also will assume that  
$\gamma$ is an oriented piecewise-smooth
curve having only transversal self-intersec\-tions and containing no poles of $P(z)$ or $q(z).$ In this paper we study the function 
$$f(t)=\sum_{i=0}^{\infty}m_it^i,$$ where 
\be \la{kuku} m_i=\int_{\gamma}P^i(z)q(z)d z.\ee
More precisely, instead of studying the function $f(t)$ directly, we study the function 
defined near infinity by the integral
\be \la{int} I(t)=
I(q,P,\gamma,t)=\frac{1}{2\pi \sqrt{-1}}\int_{\gamma} \frac{q(z)dz}{P(z)-t}.
\ee
Clearly, integral \eqref{int} defines a holomorphic function in each domain of the complement of $P(\gamma)$ in  
$\C\P^1$ and, if $I_{\infty}(t)$ is a function defined in the domain $U_{\infty}$
containing infinity, then the calculation of its Taylor series 
at infinity shows that
\be \la{prois} I_{\infty}(t)=-\frac{1}{2\pi \sqrt{-1} }\frac{1}{t}f\left(\frac{1}{t}\right).\ee Therefore, the study of $f(t)$ near zero is equivalent to the study of $I_{\infty}(t)$ near infinity and vice versa. 

Notice that under certain conditions the function 
$I_{\infty}(t)=I_{\infty}(q,P,\gamma,t)$ coincides with a similar function $I_{\infty}(\tilde q,\tilde P,\tilde \gamma,t),$ where 
$\tilde P(z),$ $\tilde q(z)$ are rational functions of smaller degrees. 
Namely, suppose that there exist rational functions $\tilde q(z),$ $\tilde P(z)$,
$W(z)$ such that $\deg W(z)>1$ and 
\be \la{poi} P(z)=\tilde P(W(z)), \ \ \ q(z)=\tilde q(W(z))W'(z).\ee
Then changing the variable $z$ to $W(z)$ we see that 
\be\la{ch1} I_{\infty}(q,P,\gamma,t)=I_{\infty}(\tilde q,\tilde P,W(\gamma),t),\ee
or equivalently
\be \la{ch2} \int _{\gamma}P^i(z)q(z)d z= \int _{W(\gamma)}\tilde P^i(z)\tilde q(z)d z, \ \ \ i\geq 0.\ee 
If rational functions $\tilde q(z),$ $\tilde P(z)$, and 
$W(z)$ as above exist, we will say that the function $I_{\infty}(q,P,\gamma,t)$ is {\it reducible} and 
{\it reduces} to the function $I_{\infty}(\tilde q,\tilde P,W(\gamma),t).$

A version of the definition above is given by the integral 
$$\frac{1}{2\pi \sqrt{-1}}\int_{\gamma} \frac{dQ(z)}{P(z)-t},$$ 
where $P(z),$ $Q(z)$ are rational functions. 
Of coarse, near infinity this integral coincides with the function 
$I_{\infty}(q,P,\gamma,t),$ where $q(z)=Q^{\prime}(z).$ Nevertheless, the definition \eqref{int} is more general since the indefinite integral $\int q(z)d z$ 
is not always a rational function. Notice that in the case where a rational function $Q(z)$ such that $q(z)=Q^{\prime}(z)$ exists, the condition \eqref{poi} is equivalent to the condition that 
\be \la{comp} P(z)=\tilde P(W(z)), \ \ \ Q(z)=\tilde Q(W(z))\ee
for some rational functions $\tilde Q(z),$ $\tilde P(z)$,
$W(z)$, $\deg W(z)>1$. If the last condition is satisfied for some rational functions $P(z),$ $Q(z)$
we will say that $P(z),$ $Q(z)$ have {\it a non-trivial common compositional right  
factor}.

Let $C_P\subset \C\P^1$ be the set of branch points of the algebraic function 
$P^{-1}(z)$ inverse to $P(z)$. Throughout the paper,  
$U$ will always denote a fixed simply connected subdomain of $\C$ such that 
$U\cap C_P=\emptyset$ 
and $\infty \in \partial U$. Notice that the condition $\infty \in \partial U$ implies that 
$U\cap U_{\infty}$ is not empty.
Since $U$ is simply connected, the condition $U\cap C_P=\emptyset$ implies that in $U$  
there exist $n=\deg P(z)$ single value analytical branches  
of $P^{-1}(z)$. We will denote these branches by $P^{-1}_i(z),$ $1\leq i \leq n.$ 
Under the analytic continuation along a closed curve  
the set $P^{-1}_i(z),$ $1\leq i \leq n,$ transforms 
to itself and this induces a homomorphism \be \l{homo} \pi_1(\C\P^1\setminus C_P,c)\rightarrow S_n,\ \ \ c\in U.\ee
The  
image $G_P$ of homomorphism \eqref{homo} is called the monodromy
group of $P(z)$. Recall that the group $G_P$ is permutation equivalent to 
the Galois group of the algebraic equation $P(x)-z=0$ over the ground field $\C(z).$

\subsection{Calculation of $I_{\infty}(t)$ for closed $\gamma$}
In this subsection we will assume that $\gamma$ is a closed curve. In this case for any point $z\in\C\P^1\setminus \gamma$ 
the winding number of $\gamma$ about $z$ is well defined. We will denote this number by
$\mu(\gamma,z).$
Let $z_1^{q}, \dots ,z_l^q$ be {\it finite} poles of $q(z).$ For $s,$ $1\leq s \leq l,$ denote
by $q_s(z)$ 
the principal part of the Laurent series of $q(z)$ at $z_s$, and 
set $$\psi_s(t)=\sum_{i=1}^n \left(\frac{q_s}{P'}\right)(P^{-1}_i(t)), \ \ \ t\in U.$$
Since $\psi_s(t)$ is invariant with respect to the action of the group 
$G_P$, this is a rational function.

Furthermore, denote by $z_1^{P}, \dots ,z_r^P$ poles of $P(z)$ in $\C\P^1$ and define $J_e,$ $1\leq e \leq r,$ as a subset of $\{1,2,\dots, n\},$ $n=\deg P(z),$ consisting of all $i\in\{1,2,\dots, n\}$ such that for $t$ close to infinity, $P^{-1}_i(t)$ is close to $z_e^P.$

\bt \la{t1} Let $\gamma$ be a closed curve. 
Then for all $t\in U\cap U_{\infty}$ close enough to infinity the equality
\begin{multline} \la{rrav}
I_{\infty}(t)=\sum_{e=1}^r \mu(\gamma,z_{e}^P)\sum_{i\in J_e}\left(\frac{q}{P'}\right)\left(P^{-1}_i(t)\right)- \sum_{s=1}^l \mu(\gamma,z_s^q)
\psi_{s}(t)+
 \\
+\frac{1}{2\pi \sqrt{-1}}\frac{1}{\left(P(\infty)-t\right)}\int_{\gamma} q(z)d z
\end{multline}
holds. 

\et 
\pr Indeed, for $t\in U\cap U_{\infty}$ 
poles of the function $$\frac{q(z)}{P(z)-t}$$ 
split into two 
groups. The first one 
contains poles of $q(z)$ while the second one contains the
points $P^{-1}_i(t),$ $1\leq i \leq n.$ Moreover, for all 
$t$ close enough to infinity these groups are disjointed 
and $P^{-1}_i(t),$ $1\leq i \leq n,$ are finite.

Clearly, for any $i,$ $1\leq i \leq n,$ we have:
$${\rm Res}_{\,P^{-1}_i(t)}\frac{q(z)}{P(z)-t}=\left(\frac{q}{P'}\right)(P^{-1}_i(t)).$$
On the other hand, for any $s,$ $1\leq s \leq l,$ we have:
$${\rm Res}_{z_s^q}\frac{q(z)}{P(z)-t}={\rm Res}_{z_s^q}\frac{q_s(z)}{P(z)-t}=
-\sum_{i=1}^n{\rm Res}_{\,P^{-1}_i(t)}\frac{q_s(z)}{P(z)-t}-{\rm Res}_{\infty}\frac{q_s(z)}{P(z)-t}=$$ $$
=-\sum_{i=1}^n\left(\frac{q_s}{P'}\right)(P^{-1}_i(t))-{\rm Res}_{\infty}\frac{q_s(z)}{P(z)-t}=-\psi_s(t)-{\rm Res}_{\infty}\frac{q_s(z)}{P(z)-t}.$$  
Furthermore, 
$${\rm Res}_{\infty}\frac{q_s(z)}{P(z)-t}=\frac{1}{P(\infty)-t}{\rm Res}_{\infty}q_s(z)=-\frac{1}{P(\infty)-t}{\rm Res}_{z_s^q}q_s(z),$$
and $$\sum_{i=1}^n\mu(\gamma,z_s^q){\rm Res}_{\infty}\frac{q_s(z)}{P(z)-t}
=-\frac{1}{P(\infty)-t}\sum_{i=1}^n\mu(\gamma,z_s^q){\rm Res}_{z_s^q}q_s(z)=$$
$$=-\frac{1}{P(\infty)-t}\sum_{i=1}^n\mu(\gamma,z_s^q){\rm Res}_{z_s^q}q(z)=
-\frac{1}{2\pi \sqrt{-1}}\frac{1}{\left(P(\infty)-t\right)}\int_{\gamma} q(z)d z.$$

Therefore,
\begin{multline}I_{\infty}(t)=\sum_{i=1}^n \mu(\gamma,P^{-1}_i(t)){\rm Res}_{\,P^{-1}_i(t)}\frac{q(z)}{P(z)-t}+\sum_{s=1}^l\mu(\gamma,z_s^q){\rm Res}_{z_s}\frac{q(z)}{P(z)-t}=\\
=\sum_{e=1}^r \sum_{i\in J_e}\mu(\gamma,P^{-1}_i(t))\left(\frac{q}{P'}\right)\left(P^{-1}_i(t)\right)- \sum_{s=1}^l \mu(\gamma,z_s^q)
\psi_{s}(t)+\\
+\frac{1}{2\pi \sqrt{-1}}\frac{1}{\left(P(\infty)-t\right)}\int_{\gamma} q(z)d z.
\end{multline}
Finally, since 
for $t$ close enough to infinity and $i\in J_e$ the equality  $\mu(\gamma,P^{-1}_i(t))=\mu(\gamma,z_e^P)$ holds we obtain \eqref{rrav}. 
\qed

Following \cite{pry}, we will say that points $x_1,x_2,\dots,x_k\in \C\P^1$ 
{\it lie on one side} of a curve $\gamma$,
if 
$\gamma$ is closed and homologous to zero in $\C\P^1\setminus \{x_1,x_2,\dots,x_k\}$. An equivalent condition 
is that  
$\gamma$ is closed and  
\be \la{hom} \mu(\gamma,x_1)=\mu(\gamma,x_2)=\dots =\mu(\gamma,x_k).\ee Further, 
if $x_1,x_2,\dots,x_k\in \C\P^1$ lie on one side of $\gamma$ and 
all numbers in \eqref{hom} equal zero, then we will say that $x_1,x_2,\dots,x_k$ {\it lie outside} $\gamma$.

\bc \la{gavv} If poles of $P(z)$ lie on one side of $\gamma$, then for any rational function $q(z)$ the function $I_{\infty}(t)$ is rational. If poles of $P(z)$ and $q(z)$ lie outside $\gamma$, then $I_{\infty}(t)\equiv 0.$
\ec
\pr Indeed, if poles of $P(z)$ lie on one side of $\gamma$, then \eqref{hom} implies that the expression  
$$\sum_{e=1}^r \mu(\gamma,z_{e}^P)\sum_{i\in J_e}\left(\frac{q}{P'}\right)\left(P^{-1}_i(t)\right)$$ in formula \eqref{rrav} is invariant with respect to the action of the group 
$G_P$ and therefore is a rational function. 
Since other terms of 
\eqref{rrav} also are rational, this implies the rationality of  
$I_{\infty}(t)$. 

The second part of the corollary follows from formula \eqref{rrav} or directly from the Cauchy theorem applied to coefficients \eqref{kuku} of $I_{\infty}(t)$.
\qed

\subsection{Calculation of $I_{\infty}(t)$ for non-closed $\gamma$} In this subsection we will assume that $\gamma$ is a non-closed 
curve with the starting point $a$ and the ending point $b$.

\bl \la{dur} The function 
$$\hat q(t)=\int_{\gamma} \frac{q(z)-q(t)}{z-t}dz
$$ is a rational function the set of poles of which
is a subset of the set of poles of $q(z).$ 
\el
\pr Indeed, $q(z)$ can be represented as a sum of terms $\alpha(z-\beta)^l,$ $l\in \Z,$ $\alpha,\beta\in \C.$ 
If $l\geq 0$, then $(z-\beta)^l-(t-\beta)^l$ is divisible by 
$z-t$ and \be \la{intt} \int_{\gamma} \frac{(z-\beta)^l-(t-\beta)^l}{z-t}dz
\ee is a polynomial.
On the other hand, 
$$(z-\beta)^{-l}-(t-\beta)^{-l}=
\frac{(t-\beta)^l-(z-\beta)^l}{(t-\beta)^l(z-\beta)^l}$$ implying that for $l<0$ integral 
\eqref{intt} is a rational function with a unique possible pole $\beta.$
\qed 

For $t\in U$ denote by 
$\Log_{1,i}(z-P^{-1}_i(t)),$ $1\leq i \leq n,$
a branch of the logarithm defined 
in a neighborhood of the point $z=a$ and by 
$\Log_{2,i}(z-P^{-1}_i(t)),$ $1\leq i \leq n,$ its analytical continuation along $\gamma$ to a neighborhood
of $b$.

\bt \la{t2} Let $\gamma$ be a non-closed curve, with the starting point $a$ and the ending point $b$,
and $t\in U\cap U_{\infty}$ be a point close enough to infinity.
Then in a neighborhood of $t$ the equality
\begin{multline} \la{rrav2}
I_{\infty}(t)=\frac{1}{2\pi \sqrt{-1}}\frac{1}{\left(P(\infty)-t\right)}\int_{\gamma} q(z)d z+\frac{1}{2\pi \sqrt{-1}}\sum_{i=1}^n \left(\frac{\hat q}{P'}\right)\left(P^{-1}_i(t)\right)+
\\ +\frac{1}{2\pi \sqrt{-1}}\sum_{i=1}^n \left(\frac{q}{P'}\right)\left(P^{-1}_i(t)\right)\left[\Log_{2,i}\left((b-P^{-1}_i(t)\right)-\Log_{1,i}\left((a-P^{-1}_i(t)\right)
\right]
\end{multline} holds.
\et 

\pr For $t$ as above 
poles $P^{-1}_i(t)$, $1\leq i \leq n,$ of 
the function $$\frac{1}{P(z)-t}$$ are finite and 
$$
{\rm Res}_{\,P^{-1}_i(t)}\frac{1}{P(z)-t}=\frac{1}{P'(P^{-1}_i(t))},\ \ \ 1\leq i \leq n.$$ This implies that  
$$\frac{1}{P(z)-t}=\frac{1}{P(\infty)-t}+\sum_{i=1}^n\frac{1}{P'(P^{-1}_i(t))}\frac{1}{\left(z-P^{-1}_i(t)\right)},$$
and hence 
$$I_{\infty}(t)=\frac{1}{2\pi \sqrt{-1}}\frac{1}{\left(P(\infty)-t\right)}\int_{\gamma} q(z)d z+\frac{1}{2\pi \sqrt{-1}}
\sum_{i=1}^n \frac{1}{P'(P^{-1}_i(t))}\int_{\gamma}\frac{q(z)d z}{z-P^{-1}_i(t)}.$$

Writing now $q(z)$ as $$q(z)=q(z)-q(P^{-1}_i(t))+q(P^{-1}_i(t)),\ \ \ 1\leq i \leq n,$$ we obtain \eqref{rrav2}. \qed

\section{Conditions for $I_{\infty}(t)$ to be rational and to vanish}
Let $P(z)$ be a rational function and $\gamma$ be a curve.
In this section we construct a finite system of equations 
\be \l{su}
\sum_{i=1}^nf_{s,i}\left(\frac{q}{P'}\right)(P^{-1}_{i}(z))=0, \ \ \ \ \ \ f_{s,i} \in \Z,  \ \ \ \ \ \ 1\leq s\leq k, 
\ee involving superpositions of the rational function $(q/P')(z)$ with branches $P^{-1}_{i}(z),$ $1\leq i \leq n,$ and 
depending on $P(z)$ and $\gamma$ only, 
such that $I_{\infty}(t)=I_{\infty}(q,P,\gamma, t)$ is rational for a rational function $q(z)$ if and only if \eqref{su} holds.
Another important result of this section 
states that if functions $P(z),$ $q(z)$ satisfy the conditions  
$q^{-1}\{\infty\}\subseteq P^{-1}\{\infty\}$ and $P(\infty)=\infty$, then the rationality of $I_{\infty}(t)$
implies that $I_{\infty}(t)\equiv 0.$

The criterion for rationality of $I_{\infty}(t)$ given below generalizes 
the criterion for orthogonality of a polynomial
$q(z)$ to all powers of a polynomial $P(z)$ given in \cite{pp}, and 
as in \cite{pp} the idea is to change the integration path $\gamma$ to a very special one. 
Let $P(z)$ be a rational function of degree $n$.
Define 
an embedded into the Riemann sphere graph $\lambda_P,$ 
associated with $P(z),$ as follows:
take a ``star'' $S$ joining a non-branch point $c$ of $P^{-1}(z)$ with all its {\it finite} branch points $c_1,c_2, ... ,c_k$ by non intersecting oriented arcs
$\gamma_1, \gamma_2, ... ,\gamma_k$, and set $\lambda_P=P^{-1}\{S\}$.
More precisely, define  
vertices of $\lambda_P$ as 
preimages
of the points $c$ and $c_s,$ $1\leq s \leq k,$ 
and edges of $\lambda_P$ as preimages of the arcs
$\gamma_s,$ $1 \leq s \leq k,$ 
under the mapping $P(z)\,:\, \C\P^1\rightarrow \C\P^1.$ 
Furthermore, for each $s,$ $1\leq s \leq k,$ mark vertices of $\lambda_P$ which are preimages of the point $c_s$ by the number $s$ (see Fig. 1). 
\begin{figure}[ht]
\epsfxsize=10.5truecm
\centerline{\epsffile{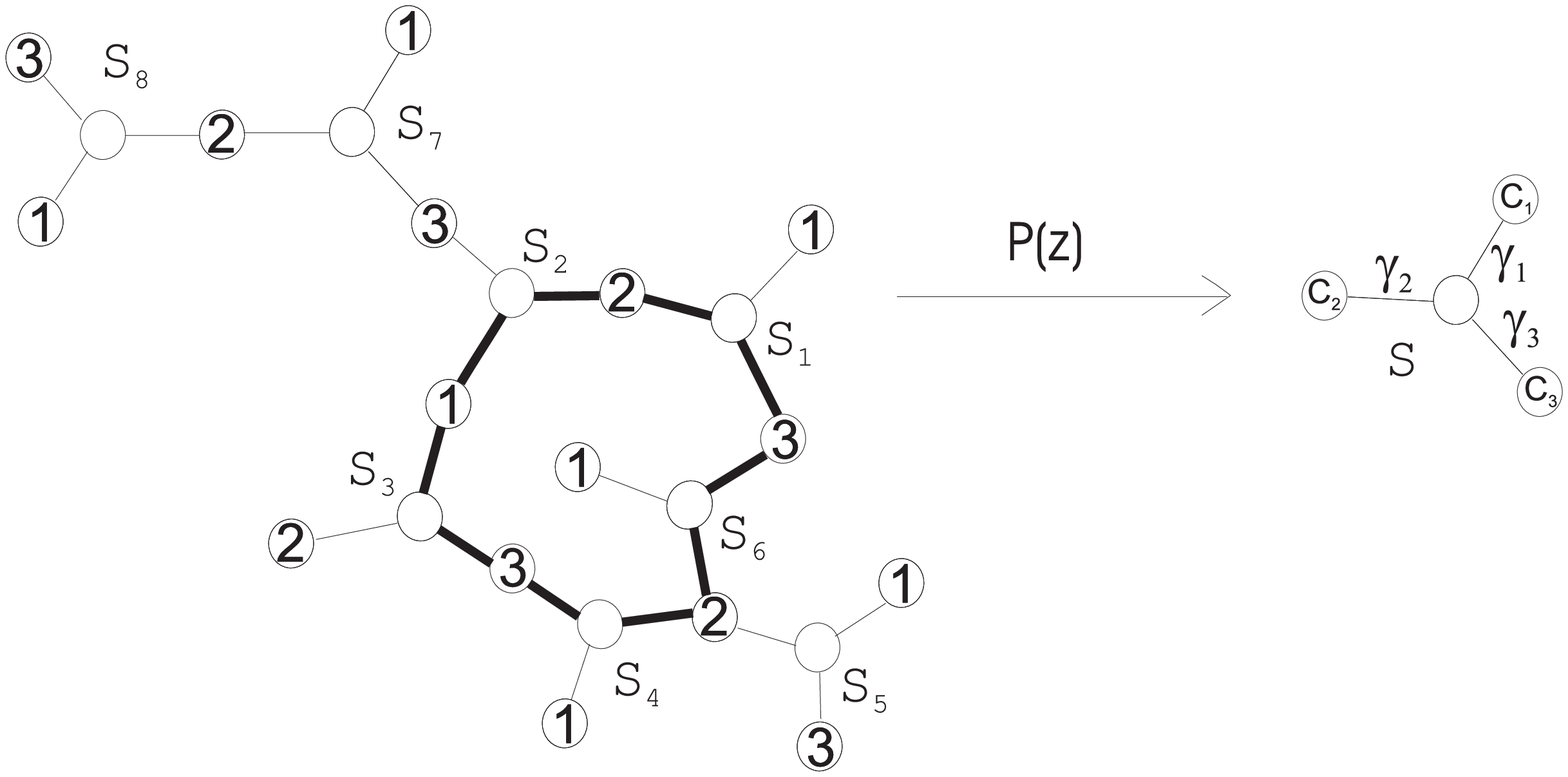}}
\smallskip
\centerline{Figure 1}
\end{figure}

By construction, the restriction of $P(z)$ on $\C\P^1\setminus \lambda_P$ is a covering of the topological punctured disk $\C\P^1\setminus \{S\cup \infty\}$ and therefore $\C\P^1\setminus \lambda_P$ is a disjointed union of punctured disks (see e.g. \cite{fors}). This implies that the graph $\lambda_P$ is connected and the faces of $\lambda_P$ are in a one-to-one correspondence with poles of $P(z).$

Define a {\it star} of $\lambda_P$ as a subset of edges of $\lambda_P$ consisting of edges adjacent to some 
non-marked vertex. Since without loss of generality we may assume 
that the domain $U$ defined above satisfies the condition  
$S\setminus\{c_1,c_2,
... ,c_k\}\subset U$, the set of stars
of $\lambda_P$ may be naturally identified with the set 
of single-valued branches of 
$P^{-1}(z)$ 
in $U$ as follows:
to the branch $P^{-1}_i(z),$ $1\leq i \leq n,$  
corresponds the star $S_i$
such that $P^{-1}_i(z)$ maps bijectively the interior of $S$ to the interior of $S_i.$

The graph constructed above is known under the name of ``constellation'' and is closely related to the notion of a 
"Dessins d'enfant" (see \cite{lz}
for further details and other versions of this construction).
Notice that the Riemann existence theorem implies that 
a rational function $P(z)$ is defined by $c_1,c_2, ... ,c_k$
and $\lambda_P$ up to a composition $P(z)\rightarrow P(\mu (z)),$ where $\mu(z)$ is a M{\"o}bius transformation.

It follows from the definition that a point $x$ is a vertex of $\lambda_P$ if and only if $P(x)$ is a branch point of $P^{-1}(z)$. 
For our purposes however it is more convenient  
to define the graph $\lambda_P$ so that in the case where the integration path $\gamma$
is non-closed, its end points $a,b$ always
would be vertices of $\lambda_P$. 
So, in the case where $\gamma$ is non-closed 
and $P(a)$ or $P(b)$ (or both of them) is not a branch point of $P^{-1}(z)$, we modify the construction as follows. Define $c_1,c_2, ... ,c_{k}$ as the set of all finite branch points of $P(z)$ supplemented by $P(a)$ or $P(b)$ (or by both of them) and set as above $\lambda_P=P^{-1}\{S\},$ where $S$ is a star
connecting $c$ with $c_1,c_2, ... ,c_{k}$.
Clearly, $\lambda_P$ is still connected and 
the points $a,b$ now are vertices of $\lambda_P.$

Deform now the integration path $\gamma$ to a path $\tilde \gamma$ 
such that $\tilde \gamma$ is contained in $\lambda_P$
and is homologous to 
$\gamma$ in $\C\P^1\setminus P^{-1}\{\infty\}$ (see Fig. 2), and 
\begin{figure}[ht]
\epsfxsize=12.2truecm
\centerline{\epsffile{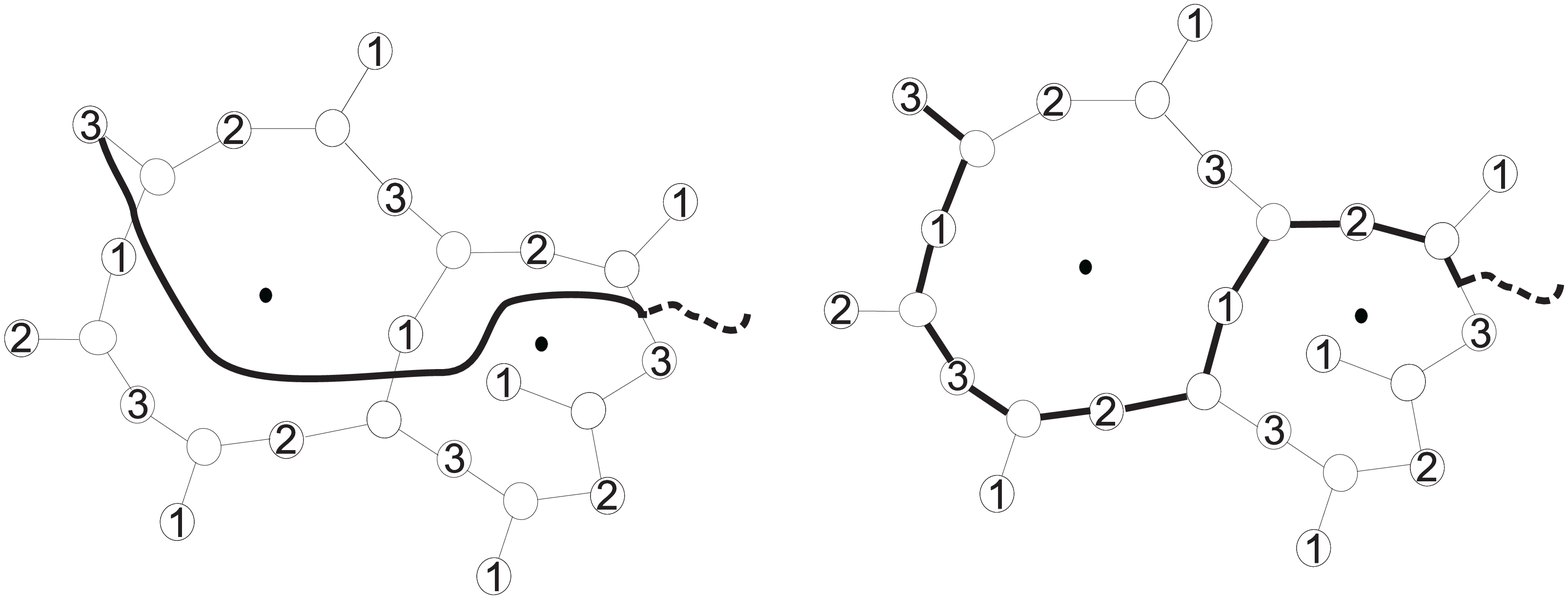}}
\smallskip
\centerline{Figure 2}
\end{figure}
define a function $\tilde I_{\infty}(t)$ by the formula 
\be \la{hh}
\tilde I_{\infty}(t)= \int_{\tilde \gamma} \frac{q(z)dz}{P(z)-t}.
\ee
Since we can write $\tilde I_{\infty}(t)$ in the form
$$\tilde I_{\infty}(t)=\int_{\tilde \gamma} \frac{q(z)d P(z)}{P^{\prime}(z)(P(z)-t)},$$
making the change of variable $z\rightarrow P(z)$
we obtain that  
\be \l{f}
\tilde I_{\infty}(t)=\sum_{s=1}^{k}\int_{\gamma_s}\frac{\phi_s(z)}{z-t}\, d z,
\ee
where $\phi_s(z),$ $1\leq s \leq k,$ are linear combinations 
of the functions $(q/P')(P^{-1}_{i}(z)),$ $1\leq i \leq n,$ in $U.$ 
Namely, for $i,$ $1\leq i \leq n,$ and $s$, $1\leq s\leq k,$ let
$c_{s,i}$ be a unique vertex of the star $S_i$ marked by the number $s.$ Then 
for $s,$ $1\leq s \leq k,$ we have: \be \l{piz}
\phi_s(z)=\sum_{i=1}^nf_{s,i}\left(\frac{q}{P'}\right)(P^{-1}_{i}(z)), 
\ee
where $f_{s,i}$ is a sum of ``signed'' appearances of 
the vertex $c_{s,i}$ 
on the path $\tilde \gamma$.
By definition, this means that an appearance
is taken with the sign plus if the center of $S_i$ is followed by $c_{s,i}$, and minus
if $c_{s,i}$
is followed by the center of $S_i$.
For example, for the graph $\lambda_P$ shown on Fig. 1 and 
the path $\tilde \gamma\subset \lambda_P$ pictured by the fat 
line we have:  
\begin{gather*}
\phi_1(z)=-\left(\frac{q}{P'}\right)(P^{-1}_{2}(z))+\left(\frac{q}{P'}\right)(P^{-1}_{3}(z)),\\
\phi_2(z)=\left(\frac{q}{P'}\right)(P^{-1}_{2}(z))-\left(\frac{q}{P'}\right)(P^{-1}_{1}(z))+\left(\frac{q}{P'}\right)(P^{-1}_{6}(z))- \left(\frac{q}{P'}\right)(P^{-1}_{4}(z)),\\
\phi_3(z)=\left(\frac{q}{P'}\right)(P^{-1}_{1}(z))-\left(\frac{q}{P'}\right)(P^{-1}_{6}(z))+ \left(\frac{q}{P'}\right)(P^{-1}_{4}(z))-\left(\frac{q}{P'}\right)(P^{-1}_{3}(z)).
\end{gather*}
\vskip 0.2cm

\bl \la{ebti} The set of equations $\phi_s(z)=0,$ $1\leq s \leq k,$
contains an equation whose coefficients are not all equal to zero 
if and only if poles $P(z)$ do not lie on one side of $\gamma$.
\el
\pr If $\gamma$ is non-closed and $x$ is the starting point or the 
ending point of $\gamma$,
then it follows from the construction that 
for $s,$ $1\leq s \leq k,$ 
such that $P(x)=c_s$ and $i,$ $1\leq i \leq n,$ such that $x\in S_i$, 
the coefficient $f_{s,i}$ of 
the equation $\phi_s(z)$ is distinct from zero. 

Assume now that $\gamma$ is a closed curve such that 
all equations $\phi_s(z)=0,$ $1\leq s \leq k,$
have zero coefficients and show that this implies that poles of $P(z)$ lie on one side of $\gamma$.
Since oriented bounds $\gamma_j,$ $1\leq j \leq d,$ of faces $f_j$ 
of the graph $\lambda_P$ generate $H_1(\C\P^1\setminus P^{-1}\{\infty\},\Z),$  
we can write $\tilde \gamma$ as a sum \be \la{repr} \tilde \gamma=\sum_i^d e_j\gamma_j,\ \ \ e_j\in \Z,\ee 
and for any $s,$ $1\leq s \leq k,$ the equation 
$\phi_s(z)$ is obtained as a sum
$$\phi_s(z)=\sum_i^d e_j\phi_{s,j}(z),$$ where $\phi_{s,j}(z)$ is an equation similar to $\phi_s(z)$
but written for $\gamma_j,$ $1\leq j \leq d.$
If in representation \eqref{repr} an index $e_{j_0}$ is distinct from zero, then for any $s,$ $1\leq s \leq k,$ the term $e_{j_0}\gamma_{j_0}$ in \eqref{repr}
gives a non-zero contribution into the equation $\phi_s(z)=0.$ Since however all coefficients of this equation are zeros, it follows from the construction that for any face $f_{i_0}$ of $\lambda_P$ adjacent to $f_{j_0}$ the equality $e_{i_0}= e_{j_0}$ holds.  
Taking into account that we can join any two faces of $\lambda_P$ by a connected chain of 
faces, this implies that all $e_j,$ $1\leq j \leq d,$ are equal. Therefore, since $\sum_i^d \gamma_j \sim 0$ in $H_1(\C\P^1\setminus P^{-1}\{\infty\},\Z)$, poles of $P(z)$ lie on one side of $\gamma$. \qed

\bt \l{t3} Let $P(z),q(z)$ be rational functions and $\gamma$ be a curve. 
Then the function $I_{\infty}(t)$ is rational if and only if 
$\phi_s(z)\equiv 0$ for all $s,$ $1\leq s \leq k.$ 
\et 

\pr Set $\hat \gamma=\gamma-\tilde \gamma$. Since $\hat \gamma$ 
is homologous to zero in $\C\P^1\setminus P^{-1}\{\infty\},$
it follows from Corollary \ref{gavv} that $\tilde I_{\infty}(t)-I_{\infty}(t)$ is a rational function, and hence  
in order to prove the theorem it is enough to prove that the function $\tilde I_{\infty}(t)$ 
is rational if and only if 
\be \l{sss}\phi_s(z)\equiv 0, \ \ \ \ 1\leq s \leq k.\ee 

It follows from Theorem \ref{t1}, Lemma \ref{dur}, and Theorem \ref{t2} that 
a multivalued analytical function $\hat I(t)$ 
obtained by the complete analytical continuation of $\tilde I_{\infty}(t)$ may ramify only at the points
$c_1, c_2, \,...\,, c_k$ or $\infty$ and that 
any other singularity of $\hat I(t)$ 
is a pole at worst. 
This implies that in order to prove the rationality of $\tilde I_{\infty}(t)$
it is enough to show that at any of points $c_1, c_2, \,...\,, c_k$ the function $\hat I(t)$ does not ramify 
and its Laurent series contains only finite number of terms with negative exponents. Indeed, if this condition is satisfied, then $\hat I(t)$ also does not ramify at $\infty$ and hence coincides in a neighborhood of $\infty$ with $I_{\infty}(t)$.
Thus, $\hat I(t)$ does not ramify in all $\C\P^1$ and its Laurent series at any point of $\C\P^1$
contains only finite number of terms with negative exponents. Therefore,
$I_{\infty}(t)$ is rational 
by the well-known characterization of rational functions.

Making if necessary a small deformation of $S\setminus \{c_1, c_2, \,...\,, c_k\}$ 
we may assume that $\lambda_P\setminus P^{-1}\{c_1, c_2, \,...\,, c_k\}$ contains no poles of $q(z)$. Assume first that the set $P^{-1}\{c_1, c_2, \,...\,, c_k\}$ also contains no poles of $q(z)$. Then representations \eqref{hh},\eqref{f} are well defined
and give an analytical extension of $\tilde I_{\infty}(t)$ to $\C\P^1\setminus S$. Keeping for this extension the same notation 
and using the well-known boundary property of
Cauchy type integrals (see e.g. \cite{mus}), we see that for any $s,$ $1\leq s \leq k,$ and any interior point $z_0$ of $\gamma_s$ the equality
\be \l{cuu}
\lim_{t \to z_0}\!\!\!\!\,^+\tilde I_{\infty}(t)-\lim_{t \to z_0}
\!\!\!\!\,^-\tilde I_{\infty}(t)=\phi_s(z_0)
\ee holds, 
where the limits are taken when $t$ approaches $z_0$ from 
the ``right''  (resp. ``left'') side of $\gamma_s$. 
Clearly, if
$\tilde I_{\infty}(t)$ is a rational function, 
then the limits above coincide for any $z_0$ and hence \eqref{sss} holds.
On the other hand, 
if \eqref{sss} holds, then
it follows directly from formula \e{f} that
$\tilde I_{\infty}(t)\equiv 0$.

If $P^{-1}\{c_1, c_2, \,...\,, c_k\}$ contains poles of $q(z),$ then 
change the path $\tilde \gamma$ in the definition of the function $\tilde I_{\infty}(t)$ as follows.
For each $s,$ $1\leq s \leq k,$ take a small loop $\delta_s$ around $c_s$ and for each $x\in P^{-1}\{c_s\}$ 
denote by $\omega_x$ the connectivity component 
of $P^{-1}\{\delta_s\}$ which bounds a domain containing $x.$
Replace now a small part of $\tilde \gamma$
near each point $x$ such that $x\in P^{-1}\{c_s\}\cap q^{-1}\{\infty\}$ 
by a part of $\omega_x$ as it shown on Fig. 3
so that the path obtained would be homologous to 
$\gamma$ in $\C\P^1\setminus P^{-1}\{\infty\}.$ 
\begin{figure}[ht]
\epsfxsize=10.5truecm
\centerline{\epsffile{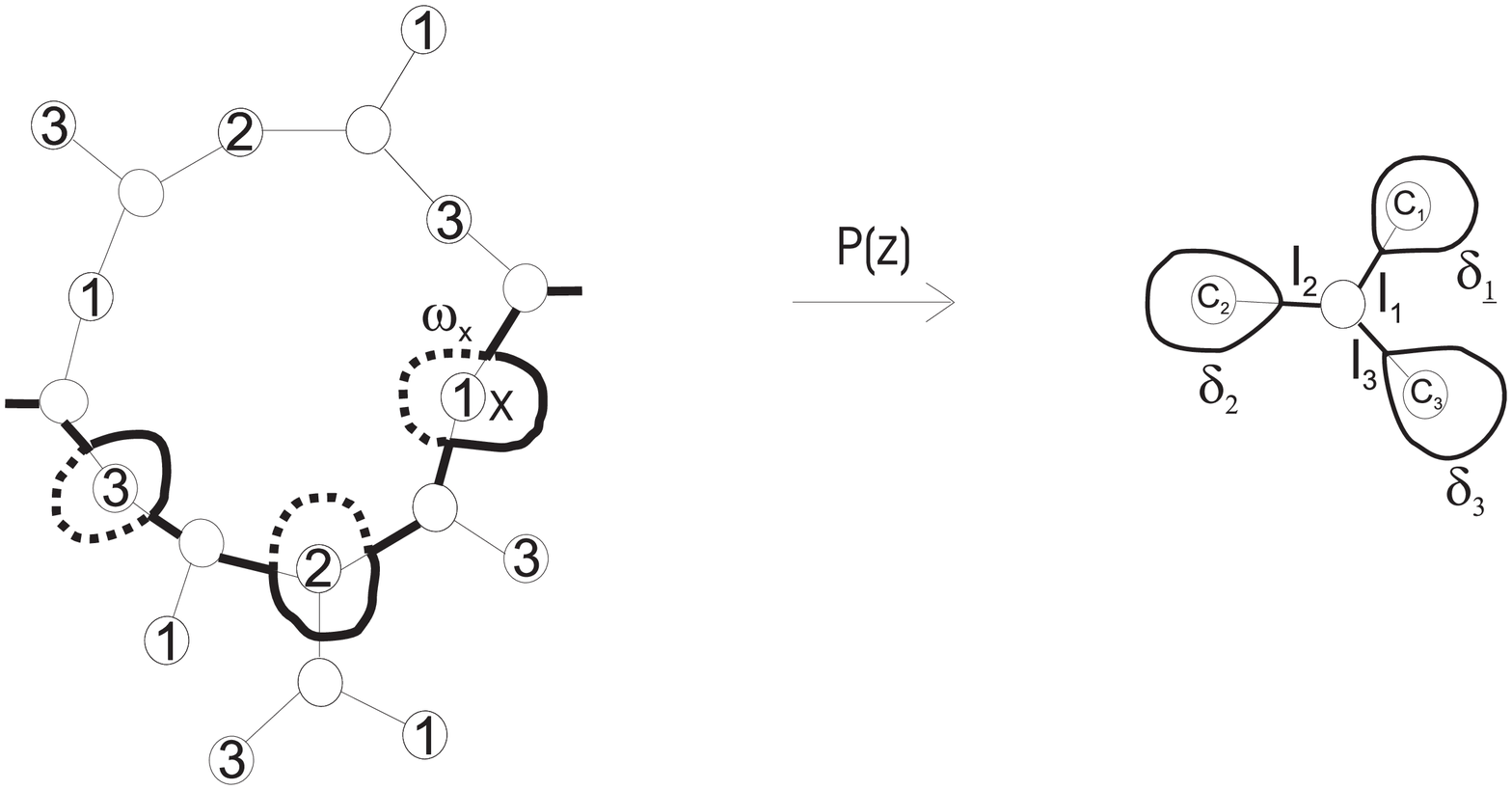}}
\medskip
\centerline{Figure 3}
\end{figure}  

Clearly, the function $\tilde I_{\infty}(t)-I_{\infty}(t)$ is still rational
and $$ \tilde I_{\infty}(t)=\sum_{s=1}^{k}\int_{l_s^+}\frac{\phi_s(z)}{z-t}\, d z
+\sum_{s=1}^{k}\int_{l_s^-}\frac{g_s(z)}{z-t}\, d z+\sum_s^{k} \int_{\delta_s}\frac{h_s(z)}{z-t}\, d z,
$$
where $l^+_s$ (resp. $l^-_s$) is the part of $\gamma_s$ which is outside (resp. inside)
the domain bounded by 
$\delta_s$, $g_s$ are some linear combinations of the functions  $(q/P')(P^{-1}_{i}(z)),$ $1\leq i \leq n,$
and $h_s(z)$ are analytic continuations of some similar combinations along $\delta_s.$ 
Since $\delta_s,$ $1\leq s \leq k,$ can be taken as close to $c_s$ as we want, formula \eqref{cuu} still implies that $\hat I(t)$ does not ramify at $c_s,$
$1\leq s \leq k,$ if and only if condition \eqref{sss} holds. In particular, condition 
\eqref{sss} is
necessary for the rationality of $\tilde I_{\infty}(t).$ 
Therefore, we only must show that if $\hat I(t)$ does not ramify in $\C\P^1$, then the Laurent series of $\hat I(t)$ at $c_s,$ $1\leq s \leq k,$ 
cannot contain an infinite number of terms with negative exponents. 

If $\gamma$ is closed, then the last statement follows directly from Theorem \ref{t1}. On the other hand, if $\gamma$ is non-closed, then it follows from Theorem 
\ref{t2} that the only points where the Laurent series of $\hat I(t)$ may have an infinite number of terms with negative exponents are points $P(a),P(b)$.
Define $J_a$ (resp. $J_b$) as a subset of $\{1,2,\dots, n\}$ consisting of all $i,$ 
$1\leq i \leq n,$ such that when $z$ is close to $P(a)$ (resp. to $P(b)$), $P^{-1}_i(z)$ is close to $a$ (resp. $b$). 

Suppose first that $P(a)\neq P(b)$. 
Then Theorem \ref{t2} implies that near $t\in U\cap U_{\infty}$
any branch of $\hat I(t)$ has the form 
$$ I_1(t)=-\frac{1}{2\pi \sqrt{-1}}\sum_{i\in J_a}\left(\frac{q}{P'}\right)\left(P^{-1}_i(t)\right)
\Log\left(a-P^{-1}_i(t)\right)+\chi(t),$$
where $\Log(z)$ is a branch of the logarithm 
and $\chi(t)$ is a branch of a function which has only a finite ramification at $P(a)$ and whose Puiseux series at $P(a)$ has only a finite number of terms with negative exponents. Furthermore, if 
$i_0\in J_a$ is a fixed index, then it is easy to see using Puiseux series that for any $i\in J_a$ we have:
$$\lim_{t \to P(a)}\frac{a-P^{-1}_i(t)}{a-P^{-1}_{i_0}(t)}\neq 0$$ and hence 
$$\Log\left(a-P^{-1}_i(t)\right)=\Log\left(a-P^{-1}_{i_0}(t)\right)+\psi_i(t),
$$ where $$\psi_i(t)=\Log\left(\frac{a-P^{-1}_i(t)}{a-P^{-1}_{i_0}(t)}\right)$$ is
a function analytical at $P(a).$ This implies that 
\be \la{kry}   I_1(t)=-\frac{1}{2\pi \sqrt{-1}}\Log\left(a-P^{-1}_{i_0}(t)\right)\sum_{i\in J_a}\left(\frac{q}{P'}\right)\left(P^{-1}_i(t)\right)
+\chi_1(t),\ee where $\chi_1(t)$ is a branch of a function which has only a finite ramification at $P(a)$ and whose Puiseux series at $P(a)$ has only a finite number of terms with negative exponents.  

Since $\Log(z)$ has an infinite ramification at $0$ 
it follows from \eqref{kry} that if $\hat I(t)$ does not ramify at $P(a)$ (or just has a finite ramification there), then necessary 
\be \la{fo1} \sum_{i\in J_a}\left(\frac{q}{P'}\right)\left(P^{-1}_i(t)\right)=0\ee implying that the Laurent series of $\tilde I_{\infty}(t)$ at $P(a)$ 
has a finite number of terms with negative exponents. 

Similarly, near $t\in U\cap U_{\infty}$
any branch of $\hat I(t)$ 
has the form 
\be \la{kry2}  I_2(t)=\frac{1}{2\pi \sqrt{-1}}\Log\left(b-P^{-1}_{j_0}(t)\right)\sum_{j\in J_b}\left(\frac{q}{P'}\right)\left(P^{-1}_i(t)\right)
+\chi_2(t),\ee 
where $j_0\in J_b$ and 
$\chi_2(t)$ is a branch of a function which has only a finite ramification at $P(b)$ and whose Puiseux series at $P(b)$ has only a finite number of terms with negative exponents. Therefore,  
if $\hat I(t)$ does not ramify at $P(b)$, then 
\be \la{fo2} \sum_{i\in J_b}\left(\frac{q}{P'}\right)\left(P^{-1}_i(t)\right)=0\ee
and the Laurent series of $\hat I(t)$ at $P(b)$ 
has a finite number of terms with negative exponents. 

Finally, if $P(a)=P(b)$, then setting $x=P(a)=P(b)$, 
we obtain similarly that 
near $t\in U\cap U_{\infty}$ any branch of $\hat I(t)$ 
has the form  
$$ I_3(t)=\frac{1}{2\pi \sqrt{-1}}\Log\left(b-P^{-1}_{j_0}(t)\right)\sum_{j\in J_b}\left(\frac{q}{P'}\right)\left(P^{-1}_i(t)\right)-$$ $$-\frac{1}{2\pi \sqrt{-1}}\Log\left(a-P^{-1}_{i_0}(t)\right)\sum_{i\in J_a}\left(\frac{q}{P'}\right)\left(P^{-1}_i(t)\right)
+\chi_3(t),$$ where $\chi_3(t)$ is a branch of a function which has only  a finite ramification at $x$  
and whose Puiseux series at $x$ has only a finite number of terms with negative exponents.  
Furthermore, if $d_a$ (resp. $d_b$) is the multiplicity of $a$ (resp. of $b$) with respect to $P(z)$, 
then for the functions $$f(t)=\left(a-P^{-1}_{i_0}(t)\right)^{d_a}, \ \ \
g(t)=\left(b-P^{-1}_{j_0}(t)\right)^{d_b}$$ the inequality $$
\lim_{t\to x}\frac{f(t)}{g(t)}\neq 0$$ holds.

Therefore, near $x$ 
$$\Log\left(b-P^{-1}_{j_0}(t)\right)=\frac{1}{d_b}\Log(g(t))=\frac{1}{d_b}\Log\left(f(t)\frac{g(t)}{f(t)}\right)=$$
$$=\frac{d_a}{d_b}\Log(a-P^{-1}_{i_0}(t))+\psi(t),
$$ where $\psi(t)=\Log(f(t)/g(t))$ is a function analytical at $x$ and hence 
\be \la{kry3}  I_3(t)=\frac{1}{2\pi \sqrt{-1}}\Log\left(a-P^{-1}_{i_0}(t)\right)\left[\sum_{j\in J_b}\left(\frac{q}{P'}\right)\left(P^{-1}_i(t)\right)-\sum_{i\in J_a}\left(\frac{q}{P'}\right)\left(P^{-1}_i(t)\right)\right]
+\ee
$$+\chi_4(t),$$ where $\chi_4(t)$ is a branch of a function which has only a finite ramification at the point $x$ 
and whose Puiseux series at $x$ has only a finite number of terms with negative exponents.
This implies that if $\hat I(t)$ does not ramify at $x$, then 
\be \la{fo3} \frac{1}{d_b}\sum_{i\in J_b}\left(\frac{q}{P'}\right)\left(P^{-1}_i(t)\right)=\frac{1}{d_a}\sum_{i\in J_a}\left(\frac{q}{P'}\right)\left(P^{-1}_i(t)\right)\ee and
the Laurent series of $\hat I(t)$ at $x$ 
has a finite number of terms with negative exponents. 
\qed

\bc \la{polusa} Suppose that poles of $P(z)$ do not lie on one side of $\gamma$. 
Then there exist integral numbers $f_i,$ $1\leq i \leq s,$ not all equal between themselves such that for any $q(z)$ such that $I_{\infty}(t)$ is rational the equality   
\be \l{pizz}
\sum_{i=1}^nf_{i}\left(\frac{q}{P'}\right)(P^{-1}_{i}(z))\equiv 0 
\ee holds.
\ec 
\pr Indeed, by construction, for any $s,$ $1\leq s \leq k,$  
either $\sum_{i=1}^nf_{s,i}=0$ or $\sum_{i=1}^nf_{s,i}=\pm 1$, where the last 
case has the place if and only if $\gamma$ is non-closed and exactly one of the end points $a$, $b$
of $\gamma$ is an $s$-vertex of $\lambda_P$. This implies that for any 
$s$ coefficients of $\phi_{s}(z)$ are not all equal, 
unless they all equal to zero. Now the lemma follows from Lemma \ref{ebti}. \qed

Notice that if \be \la{ppooll} q^{-1}\{\infty\}\subseteq P^{-1}\{\infty\},\ee 
then the condition that the curve $\tilde \gamma- \gamma$ is homologous to zero 
in $\C\P^1\setminus\{P^{-1}\{\infty\}\}$ implies that 
poles of $P(z)$ and $q(z)$ lie on one side of 
$\tilde \gamma- \gamma$. Furthermore, if additionally 
the equality 
\be \la{sku} P(\infty)=\infty\ee holds, 
then 
poles of $P(z)$ and $q(z)$ lie outside $\tilde \gamma- \gamma$ and therefore by Corollary \ref{gavv} 
the function
$\tilde I_{\infty}(t)$,
defined by
formulas \e{hh}, \e{f}, coincides with $I_{\infty}(t)$. On the other hand, condition \eqref{sss} implies that $\tilde I_{\infty}(t)\equiv 0$. Therefore, we obtain  
the following statement.

\bt \la{eg} Let $P(z),q(z)$ be rational functions and $\gamma$ be a curve. Furthermore, suppose that 
$q^{-1}\{\infty\}\subseteq P^{-1}\{\infty\}$ and $P(\infty)=\infty$.
Then the function $I_{\infty}(t)$ is rational if and only if $I_{\infty}(t)\equiv 0.$ \qed
\et

Notice that Theorem \ref{eg} also can be deduced directly from 
Theorem \ref{t1}, Lemma \ref{dur}, and Theorem \ref{t2}. Indeed, taking into account the equality $I_{\infty}(\infty)=0$, they imply that if $I_{\infty}(t)$ is rational and equalities \eqref{ppooll},
\eqref{sku} hold, then $I_{\infty}(t)$ may have poles only at points 
$c_s,$ $1\leq s \leq k.$ 

Furthermore, it follows from equalities \eqref{ppooll} and 
$$\left(\frac{1}{P'}\right)(P^{-1}_{i}(z))=(P^{-1}_{i})'(z),\ \ \ 1\leq i \leq n,$$
that 
near $c_s,$ $1\leq s \leq k,$ we have:
$$\left(\frac{q}{P'}\right)(P^{-1}_{i}(t))=O(t-c_s)^{\frac{1}{d}-1},$$ 
where $d\geq1 $ is the multiplicity of $\lim_{t\to c_s}P^{-1}_{i}(z)$ with respect to $P(z).$ 
If $\gamma$ is closed, then by Theorem \ref{t1} this implies that 
$$I_{\infty}(t)=o(t-c_s)^{-1}$$ 
near $c_s,$ $1\leq s \leq k,$ and hence if $I_{\infty}(t)$ does not ramify at $c_s,$ $1\leq s \leq k,$ then it is analytical there. Therefore, if $I_{\infty}(t)$ is rational, then it
is a constant 
equal to $I_{\infty}(\infty)=0.$ 
For non-closed $\gamma$
the proof is similar with the additional use of formulas \eqref{kry}-\eqref{fo3}.

\bc \la{kak} Let $P(z),q(z)$ be rational functions and $\gamma$ be a curve
such that $q^{-1}\{\infty\}\subseteq P^{-1}\{\infty\}$, $P(\infty)=\infty$, and   
\be \la{rave} \int_{\gamma} P^i(z) q(z) dz=0\ee 
for all $i\geq i_0,$ where $i_0\geq 0$. Then \eqref{rave} holds 
for all $i\geq 0.$ 
\ec

\pr Since \eqref{rave} implies that $I_{\infty}(t)$ is a polynomial in $1/t$, the statement follows from Theorem \ref{eg}. \qed

Notice that conditions \eqref{ppooll}, \eqref{sku}
are satisfied in particular if $P(z),$ $q(z)$ are polynomials or if 
$P(z),q(z)$ are Laurent polynomials such that $P(z)$ is not a polynomial in $z$ or in $1/z$. 

Notice also that in the course of the proof of Theorem \ref{t3} was established the 
following statement which is a version of previous results proved in \cite{pa2}, \cite{pry},
\cite{pp}.

\bp \la{egg} If $\gamma$ is non-closed and $P(a)\neq P(b),$ 
then the rationality of $I_{\infty}(t)$ implies equalities 
\eqref{fo1}, \eqref{fo2}. On the other hand, if $\gamma$ is non-closed and $P(a)=P(b),$ then 
the rationality of $I_{\infty}(t)$ implies equality \eqref{fo3}. \qed
\ep
\bc \la{rysa} If $\gamma$ is non-closed, then the rationality of $I_{\infty}(t)$ implies
that either $P(a)=P(b)$ or both points $a,b$ are ramification points of $P.$
\ec
\pr Indeed, if $P(a)\neq P(b)$ and say $a$ is not a ramification point of $P$, then the sum in \eqref{fo1} contains a single element implying $q\equiv 0.$ \qed

In conclusion of this section observe that if there exists 
a rational function $Q(z)$ such that $Q'(z)=q(z)$, then conditions \eqref{sss} implying the rationality 
of $I_{\infty}(t)$ can be written 
in a bit different form which is often used in previous publications.

\bp \la{prime} Suppose that $\gamma$ is closed and $q(z)=Q'(z)$
for some rational function $Q(z)$. Then $I_{\infty}(t)$ is rational if and only if the equalities
\be \l{piz2}
\sum_{i=1}^nf_{s,i}Q(P^{-1}_{i}(z))=0,\ \ \ 1\leq s \leq k,
\ee
where  
$f_{s,i},$ $1\leq s \leq k,$ $1\leq i \leq n,$ 
are coefficients from \eqref{piz}, hold for any choice of $Q(z)=\int q(z)d z$.

\ep 
\pr Since equations \eqref{sss} are obtained from \eqref{piz2} by derivation, condition \eqref{piz2} is clearly sufficient. Furthermore, Theorem \ref{t3} implies that if $I_{\infty}(t)$ is rational, then for any $s,$ $1\leq s \leq k,$ we have: \be \la{nez} \sum_{i=1}^nf_{s,i}Q(P^{-1}_{i}(z))=d_s,\ \ \ d_s\in \C.
\ee  On the other hand, since $\gamma$ is closed, it follows from the construction of system \eqref{sss} that 
the limit of the left
side of \eqref{nez} as $z$ tends to $c_s$ is zero. Therefore, $d_s=0$. \qed 

\bp \la{prime1} The conclusion of Proposition \ref{prime} holds also for non-closed $\gamma$ if 
$q^{-1}\{\infty\}\subseteq P^{-1}\{\infty\},$ $P(\infty)=\infty$, and $Q(z)=\int q(z)d z$ is chosen in such a way that $Q(a)=0.$
\ep

\pr Indeed, if $\gamma$ is non-closed, then calculating the limit of the left
side of \eqref{nez} as $z$ tends to $c_s$, we obtain one of the following equalities: $d_s=0$ if $c_s\neq P(a), P(b)$, 
$d_s=-Q(a)$ if $c_s=P(a)$ and $c_s\neq P(b),$  
$d_s=Q(b)$ if $c_s=P(b)$ and $c_s\neq P(a),$ or $d_s=Q(a)-Q(b)$ if $c_s=P(b)=P(a).$ Furthermore, if $I_{\infty}(t)$ is rational, then 
Theorem \ref{eg} implies that $Q(a)=Q(b)$ by \eqref{rave} taken for $i=0.$ Therefore, $d_s=0$ by the choice of 
$Q(z)=\int q(z)d z$. \qed

\section{\la{4} 
Case of generic position}
In this section we give a criterion for $I_{\infty}(t)$ to be 
rational or to vanish identically under condition that $P(z)$ is in generic position.
We start from recalling  
the following simple fact (see e.g. \cite{pm}, Lemma 2.3)
explaining the theoretic-functional meaning of 
the equality  
\be \l{compp} Q(P^{-1}_{i_1}(z))=Q(P^{-1}_{i_2}(z)),\ \ \ 1\leq i_1,i_2\leq n,\ee
where $P(z),$ $Q(z)$ are rational functions, $\deg P(z)=n.$

\bl \l{compos} Let $P(z),Q(z)$ be non-constant rational functions. Then $P(z)$ and $Q(z)$ have a non-trivial compositional right factor if and only if equality \eqref{compp} holds for some $i_1\neq i_2.$
In particular, $Q(z)=\tilde Q(P(z))$ for some rational function $\tilde Q(z)$ if and only if 
all the functions
$Q(P^{-1}_{i}(z)),$ $1\leq i \leq n,$ are equal. \qed
\el

Let $G\subset S_n$ be a transitive permutation group.  
Recall that the permutation matrix 
representation of $G$ over a field $K$ of characteristic zero is the homomorphism $\rho_G: G\rightarrow GL(K^n),$ where $\rho_G(g),$ $g\in G,$ is defined as a linear map which sends
a vector $\vec a=(a_1,a_2, ... , a_n)$ to the vector
$\vec{a_g}=(a_{g(1)},a_{g(2)}, ... , a_{g(n)}).$ Notice that for any $K$ the linear space $K^n$ has at least two $\rho_G$-invariant subspaces: 
the subspace $E_K\subset K^n$ generated by the vector
$(1,1, ... ,1)$,
and its orthogonal complement
$E_K^{\perp}$ with respect to the 
inner product 
$$(\vec u, \vec v)=u_1v_1+u_2v_2+\dots + u_nv_n,\ \ \vec v=(v_1,v_2,\dots,v_n), \ \ \vec u=(u_1,u_2,\dots,u_n).$$

\bt \la{gen} Let $P(z),$ $q(z)$ be rational functions and 
$\gamma\subset \C$ be a curve. Assume that 
$E_{\Q}$ and $E_{\Q}^{\perp}$
are the only invariant subspaces with respect to the permutation matrix representation of the monodromy group $G_P$ of $P(z)$ over $\Q.$ 
Then the function $I_{\infty}(t)$ is rational if and only if
either $\gamma$ is closed and poles of $P(z)$ lie on one 
side of $\gamma$, or $P(\gamma)$ is closed and 
$q(z)=\tilde q(P(z))P'(z)$ 
for some rational function $\tilde q(z)$. 
\et
\pr
The sufficiency follows from Corollary \ref{gavv} taking into account that if 
$q(z)=\tilde q(P(z))P'(z)$, then  
\be \la{uzx} I_{\infty}(q,P,\gamma,t)=I_{\infty}(\tilde q,z,P(\gamma),t).\ee

Assume now that $I_{\infty}(t)$ is rational and 
poles of $P(z)$ do not lie on one side of $\gamma$. In this case by Corollary \ref{polusa} 
there exist integers $f_i,$ $1\leq i \leq s,$ not all equal between themselves such that 
equality \eqref{pizz} holds.
Furthermore, continuing analytically 
equality \eqref{pizz} and exchanging $\sigma$ and $\sigma^{-1}$, we see that
for any $\sigma\in G_P$ the equality 
$$ \sum_{i=1}^{n} f_{\sigma(i)}\left(\frac{q}{P'}\right)(P_i^{-1}(z))=0$$ holds.

Let $V$ be a subspace
of $\Q^n$ generated by the vectors $\vec v_{\sigma},$
$\sigma\in G_P,$ where $$\vec v_{\sigma}=(f_{\sigma(1)},
f_{\sigma(2)}, ...
, f_{\sigma(n)}).$$ Clearly, $V$ is $\rho_{G_P}$-invariant and for any $\vec v=(v_1, v_2, \dots, v_n)$ from $V$ 
the equality 
\be \la{pizz1} \sum_{i=1}^{n} v_i\left(\frac{q}{P'}\right)(P_i^{-1}(z))=0\ee holds.
Since $f_i,$ $1\leq i \leq s,$ are not all equal between themselves, $V\neq E_{\Q}$. Furthermore, $V\neq \Q^n$. Indeed, otherwise 
the elements $\vec {e}_i,$ $1\leq i \leq n,$ of the Euclidean basis are contained in $V$, and \eqref{pizz1} implies that  
$$\left(\frac{q}{P'}\right)(P_{1}^{-1}(z))\equiv \left(\frac{q}{P'}\right)(P_{2}^{-1}(z))\equiv \dots \equiv \left(\frac{q}{P'}\right)(P_{n}^{-1}(z)\equiv 0$$ in contradiction with the assumption $q(z)\not\equiv 0$. Hence, $V=E^{\perp}_{\Q}.$ 
Since this implies that $V$ contains the vectors $\vec {e}_i-\vec {e}_j,$ $1\leq i,j \leq n,$ it follows from \eqref{pizz1} that

\be \la{poio} \left(\frac{q}{P'}\right)(P_{1}^{-1}(z))\equiv \left(\frac{q}{P'}\right)(P_{2}^{-1}(z))\equiv \dots \equiv \left(\frac{q}{P'}\right)(P_{n}^{-1}(z)).\ee
Therefore, 
$\left(\frac{q}{P'}\right)(z)=\tilde q(P(z))$ for some rational function $\tilde q(z)$ by Lemma \ref{compos}. 
Finally, 
it follows from equality \eqref{uzx} and Corollary \ref{rysa} 
that $P(\gamma)$ is closed.
\qed

Recall that a permutation group $G$ acting on a set $C$ is called doubly transitive if it acts transitively on the set of pairs of elements of $C$. 
Notice that the full symmetric group is obviously doubly transitive.

\bc \la{genn} Let $P(z),$ $q(z)$ be rational functions and 
$\gamma\subset \C$ be a curve. Assume that the monodromy group $G_P$ of $P(z)$ is doubly transitive.
Then the function $I_{\infty}(t)$ is rational if and only if
either $\gamma$ is closed and poles of $P(z)$ lie on one 
side of $\gamma$, or $P(\gamma)$ is closed and 
$q(z)=\tilde q(P(z))P'(z)$ 
for some rational function $\tilde q(z)$.
\ec
\pr Indeed, it is well known
(see e.g. \cite{wi}, Th. 29.9) that a permutation group $G$ is doubly transitive if and 
only if the subspaces $E_{\C}$ and $E_{\C}^{\perp}$ are the only $\rho_G$-invariant subspaces
with respect to the permutation matrix representation of $G_P$ over $\C.$ 
Therefore, if $G_P$ is doubly transitive, then  
$E_{\Q}$ and $E_{\Q}^{\perp}$
are the only invariant subspaces with respect to the permutation matrix representation of $G_P$ over $\Q$ and the corollary follows from Theorem \ref{gen}. 
\qed

\bc \la{gen+} There exists a proper algebraic subset 
$\Sigma\subset \C\P^{2n+1}$ such that for any rational function \be \la{rraatt} P(z)=\frac{a_nz^n+a_{n-1}z^{n-1}+...+a_0}{b_nz^n+b_{n-1}z^{n-1}+...+b_0}\ee
with $(a_n,...,a_0,b_n, ... , b_0)\notin \Sigma$, any
non-zero rational function $q(z)$, and any curve $\gamma\subset \C$, the function $I_{\infty}(t)$ is rational if and only if either $\gamma$ is closed and 
poles of $P(z)$ lie on one 
side of $\gamma$, or $P(\gamma)$ is closed and
$q(z)=\tilde q(P(z))P'(z)$ for some rational function $\tilde q(z)$.
\ec
\pr It is easy to see that there exists a proper algebraic subset 
$\Sigma\subset \C\P^{2n+1}$ 
such that corresponding rational functions
are of degree $n$, indecomposable, and have only simple branch points (a branch point $x$ 
of a rational function $f(z)$ of degree $n$ is called simple if $f^{-1}\{x\}$ 
contains $n-1$ points). This means that monodromy groups of these functions 
are primitive and contain a transposition. 
Since a primitive permutation group containing 
a transposition is the full symmetric group
(see e.g. Theorem 13.3 of \cite{wi})
the corollary follows now from Corollary \ref{genn}. \qed

\bt \la{iioo} Let $P(z),q(z)$ be rational functions and $\gamma\subset \C$ be a curve such that 
\be \la{rav}
\int_{\gamma} P^i(z) q(z) dz=0,  \ \ \ i\geq 0.
\ee 
Suppose additionally that poles of $P(z)$ do not lie on one side of $\gamma$ and that $E_{\Q}$ and $E_{\Q}^{\perp}$
are the only invariant subspaces with respect to the permutation matrix representation of the monodromy group $G_P$ of $P(z)$ over $\Q.$
Then $P(\gamma)$ is closed and $q(z)=\tilde q(P(z))P'(z)$ for some rational function $\tilde q(z)$ whose poles lie outside the curve $P(\gamma).$
\et
\pr Indeed, by Theorem \ref{gen} the curve $P(\gamma)$ is closed and $q(z)=\tilde q(P(z))P'(z)$ for some rational function $\tilde q(z)$. Furthermore, it follows from \eqref{uzx} and Theorem \ref{t1} that 
$$I_{\infty}(t)=-\sum_{s=1}^{\tilde l} \mu(P(\gamma),z_s^{\tilde q})
\tilde q_s,$$ 
where $z_1^{\tilde q},z_2^{\tilde q}, \dots, z_{\tilde l}^{\tilde q}$ are finite poles of $\tilde q(z)$ and $\tilde q_s(z)$ is  
the principal part of $\tilde q(z)$ at $\tilde z_s$. 
Therefore, the equality 
$I_{\infty}(t)=0$ implies that $\mu(P(\gamma),z_s^{\tilde q})=0,$ $1\leq s \leq \tilde l.$ \qed

The Corollaries \ref{iiooo} and \ref{gennn} below are obtained from Theorem \ref{iioo} in the same way as Corollaries \ref{genn} and \ref{gen+} are obtained from Theorem \ref{gen}.

\bc \la{iiooo} Let $P(z),q(z)$ be rational functions and $\gamma\subset \C$ be a curve such that equalities \eqref{rav}
hold.
Suppose additionally that poles of $P(z)$ do not lie on one side of $\gamma$ and that the monodromy group $G_P$ of $P(z)$ is doubly transitive.
Then $P(\gamma)$ is closed and $q(z)=\tilde q(P(z))P'(z)$ for some rational function $\tilde q(z)$ 
whose poles lie outside the curve $P(\gamma).$ \qed
\ec

\bc \la{gennn} There exists a proper algebraic subset 
$\Sigma\subset \C\P^{2n+1}$ such that for any non-zero rational function $q(z)$, any curve $\gamma\subset \C$, and any 
rational function $P(z)$ the poles of which do not lie on one side of $\gamma$ and $(a_n,...,a_0,b_n, ... , b_0)\notin \Sigma$,  
the equalities \eqref{rav} 
imply that $P(\gamma)$ is closed and $q(z)=\tilde q(P(z))P'(z)$ for some rational function $\tilde q(z)$ whose poles lie outside the curve $P(\gamma).$ \qed
\ec

\noindent{\bf Remark.} Notice that if $P(z)$ is a {\it polynomial}, then the
requirement of Theorem \ref{gen} imposed on $\rho_{G_P}$-invariant subspaces of $\Q^n$ may be weakened to the requirement of indecomposability of $P(z)$ via the Schur theorem (see \cite{pa2}). However, there exist indecomposable rational functions $P(z)$
for which Theorem \ref{gen} fails to be true (see Section \ref{cri}). 


\section{Double moments of rational functions}
In this section 
we prove two results which can be considered as versions of
the Wermer theorem \cite{w1}, \cite{w2}, describing 
analytic functions on $S^1$ 
satisfying $$\int_{S^1}h^i(z)g^j(z)g'(z)d z=0, \ \ \ i,j\geq 0,$$ in the case where the functions $h(z),g(z)$ are rational
while the integration path is allowed to be an arbitrary curve in 
$\C.$
These results also generalize Theorem 6.1 and Corollary 6.2 of \cite{pry}.

For given $P(z),$ $Q(z)$ and $j\geq 1$ denote by $I_j(t)$ the generating functions for the sequence of the moments 
\be \la{mome2}  
m_i=\int_{\gamma}P^i(z)Q^j(z)Q'(z)\,dz, \ \ \ i\geq 0.
\ee
\bt \la{wermer} Let $P(z),Q(z)$ be rational functions and 
$\gamma$ be a curve
such that the functions $I_j(t)$ are rational for any 
$j,$ $j_0\leq j \leq j_0+ n-1,$ where $j_0\geq 0$ and $n=\deg P(z).$
Then there exist rational functions $\tilde P(z),$
$\tilde Q(z),$ $W(z)$ such that: 
\be \la{comp1} P(z)=\tilde P(W(z)),\ \ \ Q(z)=\tilde Q(W(z)),\ee the curve $W(\gamma)$ is 
closed, and poles of $\tilde P(z)$ lie on one side of $W(\gamma)$.
\et

\vskip 0.2cm

\pr 
The proof uses the same idea as the proof of Theorem 2 in \cite{pa2}, where double moments
of polynomials were investigated. 
Let $W(z)$ be a rational function such that $\C(P(z),Q(z))=\C(W(z)).$ 
Then 
the corresponding functions $\tilde P(z),$ $\tilde Q(z)$ in \eqref{comp1} have no common compositional right factor. 
Since for any $i\geq 0,j\geq 0$ we have:
$$
\int_{\gamma} P^i(z)Q^j(z)Q'(z)\,dz=\int_{W(\gamma)} \tilde P^i(z)\tilde Q^j(z)\tilde Q'(z)\,dz,
$$
it is enough to show that $\tilde \gamma=W(\gamma)$ is closed and 
poles of $\tilde P(z)$ lie on one side of $\tilde \gamma.$ Assume 
the inverse. Then it follows from Corollary \ref{polusa} 
applied to $\tilde P(z)$ and 
$\tilde Q^j(z)\tilde Q^{\prime}(z),$ $j_0 \leq j \leq j_0+\tilde n-1,$
where $\tilde n=\deg \tilde P(z)\leq n,$ 
that the system 
\be \la{sys} \sum_{s=1}^{\tilde n} \tilde f_i\tilde Q^j(\tilde P_{i}^{-1}(z))\left(\frac{\tilde Q^{\prime}}{P'}\right)(\tilde P_{i}^{-1}(z))=0, \ \ \ j_0 \leq j \leq j_0+\tilde n-1,\ee
has a non-trivial solution $\tilde f_1,\tilde f_2,\dots \tilde f_{\tilde n}$. 
Since the determinant of \eqref{sys} is a product of the 
Vandermonde determinant $D=\parallel \tilde Q^j(\tilde P^{-1}_i(z))\parallel$ and a non-zero function 
$$\prod_{i=1}^{\tilde n}\tilde Q^{j_0}(\tilde P^{-1}_{i}(z)) \left(\frac{\tilde Q^{\prime}}{P'}\right)(\tilde P_{i}^{-1}(z)),$$ this implies that
\be \la{comp2} \tilde Q(\tilde P^{-1}_{i_1}(z)) \equiv \tilde Q(\tilde P^{-1}_{i_2}(z))\ee for some $i_1\neq i_2,$ $1\leq i_1,i_2 \leq \tilde n,$ and hence 
$\tilde P(z),$ $\tilde Q(z)$ have a common compositional right factor
by Lemma \ref{compos}. The obtained contradiction proves the theorem. \qed

\bt Let $P(z),Q(z)$ be rational functions and 
$\gamma$ be a curve. Then equalities 
\be \la{wer4} \int_{\gamma}
P^i(z)Q^{j}(z)Q^{\prime}(z)d z=0 \ee
hold for all $i\geq i_0,$ $j\geq j_0,$
where $i_0\geq 0,$ $j_0\geq 0,$  
if and only if there exist rational functions $\tilde P(z),$
$\tilde Q(z),$ $W(z)$ such that: 
\be \la{comp3} P(z)=\tilde P(W(z)),\ \ \ Q(z)=\tilde Q(W(z)),\ee the curve $W(\gamma)$ is 
closed, and poles of $\tilde P(z)$ and $\tilde Q(z)$ lie on one side of the curve $W(\gamma)$. In particular, if equalities \eqref{wer4} hold for all $i\geq i_0,$ $j\geq j_0,$ then they hold for all $i\geq 0,$ $j\geq 0.$ 
\et

\pr Using the same reduction as in the proof of Theorem \ref{wermer}, 
it is enough to show that if $P(z)$ and $Q(z)$ have no common compositional right factor,
then \eqref{wer4} holds if and only if poles of $P(z)$ and $Q(z)$ lie on one side of $\gamma$. 

Assume first that poles of $P(z)$ and $Q(z)$ lie on one side of $\gamma$ 
and show that this implies that $I_j(t)=0$ for any $j\geq 0.$
Set $$q(z)=Q^j(z)Q^{\prime}(z), \ \ \ q_{\infty}(z)=q(z)-\sum_s^lq_s(z),$$ 
where $q_s(z),$ $1\leq s \leq l,$ are  
principal parts of $q(z)$ at its finite poles $z_s,$ $1\leq s \leq l.$ 
Applying Theorem \ref{t1}
to $q(z)$ 
and taking into account that the equality \be \la{rtygf} q(z)=\left(\frac{Q^{j+1}(z)}{j+1}\right)^{\prime}\ee implies 
the equality $\int_{\gamma} q(z)d z=0,$ we see 
that \be \la{blev} I_j(t)=\mu\sum_{i=1}^n\left(\frac{q_{\infty}}{P'}\right)(P^{-1}_i(t)), \ee
where $\mu$ equals to the common
winding number of poles of $P(z)$ and $Q(z).$
If $\infty$ is a pole of $Q(z)$, then $\mu=0$ implying 
$I_j(t)=0.$ On the other hand, if $\infty$ is not a pole of $Q(z)$,
then \eqref{rtygf} implies that $q_{\infty}(z)=0$.

In other direction, assume that \eqref{wer4} holds and show that 
then points from the set $P^{-1}\{\infty\}\cup Q^{-1}\{\infty\}$ lie on one side of $\gamma.$
Clearly, for any $s\geq i_0$ the function $$R(z)=P^s(z)+Q(z)$$
satisfies the equalities
$$ \int_{\gamma}
R(z)^iQ^j(z)Q^{\prime}(z)d z=0, \ \ i\geq 0, \ j\geq j_0 .$$
Furthermore, $$R^{-1}\{\infty\}\subseteq P^{-1}\{\infty\}\cup Q^{-1}\{\infty\}$$ and, 
if $s$ is big enough, then $$R^{-1}\{\infty\}= P^{-1}\{\infty\}\cup Q^{-1}\{\infty\}.$$ Therefore, since we may apply Theorem \ref{wermer} to the functions $R(z),$ $Q(z)$,
it is enough to prove that for any $k\geq 1$ 
there exists $s\geq k$ such that 
$R(z)$ and $Q(z)$, or equivalently $P^s(z)$ and $Q(z)$, have no common compositional right factor. 

Since the monodromy group of a rational function has only 
finite number of imprimitivity systems, there exist a finite number of right factors 
$Q_j(z),$ $\deg Q_j(z)>1,$ of $Q(z)$ such that any other right factor of $Q(z)$ of degree greater than one has the form $\mu\circ Q_j(z)$ for 
some $Q_j(z)$ and a M\"obius transformation $\mu.$ Further, it is easy to see that 
if $s_j$ is a minimal number such that $P^{s_j}(z)$ and $Q(z)$ have a common compositional right factor $\mu\circ Q_j(z)$, then any other $s$ with such a property is a multiple of 
$s_j.$ Since $P(z)$ and $Q(z)$ have no common compositional right factor, all $s_j$ are 
greater than one.
Therefore, for $s$ which is not a multiple of any $s_j$ the polynomials
$R(z)$ and $P(z)$ have no common compositional right factor. 
\qed

\section{Laurent polynomial moment problem}  
In this subsection we study the following problem:
{\it for a given Laurent polynomial $L(z)$ 
describe Laurent polynomials $m(z)$ such that}
\be\la{1l}
\int_{S^1} L^i(z)m(z) d z= 0,
\ee 
{\it for all $i\geq i_0,$ where $i_0\geq 0$}. 
It is easy to see that if $L(z)$ is a polynomial in $z$ or in $1/z$, then for any $m(z)$ there exists $i_0\geq 0$ such that \eqref{1l} holds.  So, the interesting case is the one where 
$L(z)$ is not a polynomial in $z$ or in $1/z$. We will call such Laurent polynomials {\it proper}.

We start from a generalizations of the
following result proved by Duistermaat and van der Kallen \cite{dui}:
if all integral positive powers of a Laurent polynomial $L(z)$ have no constant term, then
$L(z)$ is either a polynomial in $z$, or a polynomial in $1/z.$
Clearly, the condition that all powers of $L(z)$ have no constant term
is equivalent to the condition that integrals in \eqref{1l} vanish for 
$m(z)=1/z$ and $i\geq 1.$

\bt \l{d1} Let $L(z)$ and
$m(z)$ be Laurent polynomials 
such that the coefficient of the term $1/z$ in $m(z)$ is distinct from zero and \eqref{1l} holds 
for all $i\geq i_0,$ where $i_0\geq 0.$ Then  
$L(z)$ is either a polynomial in $z$, or a polynomial in $1/z.$\et
\pr Assume the inverse. Then by Corollary \ref{kak}
equalities \eqref{1l}
hold for all $i\geq 0.$ On the other hand, for $i=0$ the integral in \eqref{1l}  
does not vanish since it
coincides with the coefficient of the term $1/z$ in $m(z)$ multiplied by $2\pi \sqrt{-1} $. 
\qed

\bc \la{d} Let $L(z)$ be a proper Laurent polynomial and $m(z)$ a Laurent polynomial 
such that \eqref{1l} holds for all $i\geq i_0,$ where $i_0\geq 0.$ Then  
there exists a Laurent polynomial $M(z)$ such that $m(z)=M'(z).$ \qed
\ec

Let 
$$L(z)=a_{n_1}z^{n_1}+a_{n_1+1}z^{n_1+1}+\dots +a_{n_2}z^{n_2},\ \ \ a_{n_1}\neq 0, \ \ \ a_{n_2} \neq 0,$$ be a Laurent polynomial.
Define the {\it bi-degree} of $L(z)$ as the ordered pair $(n_1,n_2)$ of integers $n_1,n_2.$ Notice that if
$M(z)$ is another Laurent polynomial whose bi-degree is $(m_1,m_2)$, then the bi-degree of the product $L(z)M(z)$ is $(n_1+m_1,n_2+m_2).$

The next result provides yet another generalization of the theorem of Duistermaat and van der Kallen.
\bt \la{d2} Let $L(z)$ be a Laurent polynomial of bi-degree $(n_1,n_2)$ and 
$m(z)$ be either a polynomial in $z$ of bi-degree $(m_1,m_2)$, 
where 
$m_1\equiv -1\, \mod n_1,$
or a polynomial in $1/z$ of bi-degree $(m_1,m_2),$
where $m_2\equiv -1\, \mod n_2,$ such that  
\eqref{1l} holds 
for all $i\geq i_0$, where $i_0\geq 0.$ Then  
$L(z)$ is either a polynomial in $z$, or a polynomial in $1/z.$
\et
\pr Assume the inverse. Observe that then in particular $n_1<0$ and $n_2>0.$ Furthermore, Corollary \ref{kak} implies that
\eqref{1l} hold for all $i\geq 0$. Therefore, in order to prove the theorem 
it is enough to show that if $m(z)$ has the form as above, then
there exists $k\geq 0$ such that for $i=k$ integral in \eqref{1l} 
is distinct from zero.

If $m(z)$ is a polynomial in $z$ and $l_1\geq 0$ is a number such $m_1+n_1l_1=-1$, then the integral in \eqref{1l} 
is distinct from zero for $i=l_1$ since 
the bi-degree of $L^{l_1}(z)m(z)$ is $(-1,m_2+n_2l_1)$ implying that the residue of $L^{l_1}(z)m(z)$ at zero does not vanish.
Similarly, if $m(z)$ is a polynomial in $1/z$ and $l_2\geq 0$ is a number such $m_2+n_2l_2=-1$, then the integral in \eqref{1l} 
is distinct from zero for $i=l_2$, since the bi-degree of $L^{l_2}(z)m(z)$ is $(m_1+n_1l_2,-1)$. \qed 

\bc Let $L(z)$ be a Laurent polynomial of bi-degree $(n_1,n_2)$ and $d$ be either a non-negative integer such that $d\equiv 0\, \mod n_2$, or a non-positive integer such that 
$d\equiv  0\, \mod n_1$. Suppose that for all $i\geq i_0$, where $i_0\geq 0,$ 
the coefficient of $z^d$ in $L^i(z)$ vanishes. Then  
$L(z)$ is either a polynomial in $z$, or a polynomial in $1/z.$ \qed

\ec

For a Laurent polynomial 
$M(z)$ denote by $M_0(z)$ the principal part of $M(z)$ at zero and 
by $M_{\infty}(z)$ the difference $M(z)-M_0(z).$
Taking into account Corollary \ref{d} in the following we usually will write system \eqref{1l} in the form 
\be\la{1l+}
\int_{S^1} L^i(z)d M(z)= 0,
\ee $i\geq i_0,$ where it is always assumed that $M_{\infty}(0)=0.$

Define $J_0$ (resp. $J_{\infty}$) as a subset of $\{1,2,\dots, r\},$ $r=\deg L(z),$ consisting of all $i\in\{1,2,\dots, r\}$ such that for $t$ close to infinity, $L^{-1}_i(t)$ is close to $0$ 
(resp. to $\infty$). Notice that $\{1,2,\dots, r\}=J_0\cup J_{\infty}.$
The theorem below summarizes general results about $I_{\infty}(t)$ 
obtained above in the particular case where $I_{\infty}(t)$
corresponds to moments in \eqref{1l+}.

\bt \la{laurent} Let $L(z)$ be a proper Laurent polynomial and $M(z)$ is a Laurent polynomial such that \eqref{1l+} holds for all $i\geq i_0$, where $i_0\geq 0.$ Then \eqref{1l+} holds for all $i\geq 0$. Furthermore, condition \eqref{1l+} is equivalent to the condition 
\be \la{lau}
\sum_{i\in J_0}M_{\infty}(L^{-1}_i(t))\equiv 
\sum_{i\in J_{\infty}}M_0(L^{-1}_i(t)).
\ee 
Finally, if the monodromy group of $L(z)$ is doubly transitive, 
or more generally, if $E_{\Q}$ and $E_{\Q}^{\perp}$
are the only invariant subspaces with respect to the permutation matrix representation of the monodromy group $G_L$ of $L(z)$ over $\Q,$ then \eqref{1l+} holds if and only if  
$M(z)=\tilde M(L(z)),$ where $\tilde M(z)$ is a polynomial. 
\et 
\pr 
The first statement follows from Corollary \ref{kak}. Furthermore, it follows from Theorem \ref{t1}, 
taking into account the equality $\int_{S_1}d M(z)=0,$ 
that 
\begin{multline}\la{eps} I_{\infty}(t)=\sum_{i\in J_0}\left(\frac{M'}{L'}\right)(L^{-1}_i(t))
-\sum_{i=1}^{\deg L}\left(\frac{M_0^{\prime}}{L'}\right)(L^{-1}_i(t))=\\
=\sum_{i\in J_0}\left(\frac{M_{\infty}^{\prime}}{L'}\right)(L^{-1}_i(t))-
\sum_{i\in J_{\infty}}\left(\frac{M_{0}^{\prime}}{L'}\right)(L^{-1}_i(t)).
\end{multline}
Integrating the last equality we see 
that condition \eqref{1l+} is equivalent to the condition  
\be \l{ggss} \sum_{i\in J_0}M_{\infty}(L^{-1}_i(t))- 
\sum_{i\in J_{\infty}}M_0(L^{-1}_i(t))=c,\ee where 
$c\in \C.$ 
Furthermore, since the limit of the left part of \eqref{ggss} as $t$ tends to infinity
is $M_{\infty}(0)\vert J_0\vert=0$, we conclude that  
\eqref{1l+} is equivalent to \eqref{eps}.

Finally, Theorem \ref{iioo} implies that if $E_{\Q}$ and $E_{\Q}^{\perp}$
are the only invariant subspaces with respect to the permutation matrix representation of $G_L$ of $L(z)$ over $\Q,$ then there 
exists a rational function $N(z)$ such that 
\be \la{swed} M'(z)=N(L(z))L'(z).\ee Since $M'(z)$ is a Laurent polynomial it follows from \eqref{swed} that $N(L(z))$ also is a Laurent polynomial. Therefore,
$N(z)$ is a polynomial for otherwise $N(L(z))$ would have a pole distinct from $0,\infty,$ and hence $M(z)=\tilde M(L(z)),$ where $\tilde M(z)=\int N(z) dz.$ 
\qed

Notice that if 
$L(z)$ is {\it decomposable}, then Laurent polynomials $M(z)$ 
satisfying \eqref{1l+} but distinct from the ones described in Theorem \ref{laurent} 
always exist. 
Indeed, observe first that if $L(z)=A(B(z))$ is a decomposition of a Laurent polynomial $L(z)$
into a composition of rational functions $A(z)$ and $B(z)$, with $\deg A(z)>1,$ $\deg B(z)>1,$  
then the condition $B^{-1}\{A^{-1}\{\infty \}\}=\{0,\infty\},$ implies that there exists a M\"obius transformation
$\mu(z)$ such that either $A(\mu (z))$ is a polynomial and $\mu^{-1}(B(z))$ is a Laurent polynomial, or $A(\mu (z))$ is a proper Laurent polynomial and $\mu^{-1}(B(z))=z^d,$ for some $d>1$ 
\footnote{For a comprehensive decomposition theory of Laurent polynomials generalizing the decompositions theory of polynomials developed by Ritt \cite{r1} we refer the reader to \cite{pak}}.
Therefore, if $L(z)$ is decomposable, then either there exist a polynomial $\tilde L(z)$ and a Laurent polynomial
$L_1(z)$ such that $L(z)=\tilde L(L_1(z)),$ or there exists a Laurent polynomial $L_1(z)$ such that $L(z)=L_1(z^d)$ for some $d>1$. 
In the first case 
it is easy to see that 
any Laurent polynomial $M(z)=\tilde M(L_1(z)),$ where $\tilde M(z)$ is a polynomial, satisfies \eqref{1l+} for all $i\geq 0.$ On the other hand, in the second case the residue calculation shows that 
any Laurent polynomial $M(z)$, containing no terms of degrees
which are multiples of $d$,
satisfies \eqref{1l+}. Furthermore, if $L(z)$ admits several decompositions, then
the sum of corresponding $M(z)$ also satisfies \eqref{1l+}. 

It seems natural to start the investigation of solutions of the Laurent polynomial moment problem from describing {\it polynomial} solutions,
and two theorems below are initial results in this direction. Another interesting 
subproblem is to describe solutions of the polynomial moment problem in the case
where $L(z)$ is indecomposable. In the last case ``expectable'' solutions should 
have the form $M(z)=\tilde M(L(z)),$ where $\tilde M(z)$ is a polynomial. However, one can show  (see Section 8 and the paper \cite{ppz})
that other solutions also may exist.

\bt \la{d3} Let $L(z)$ be a proper Laurent polynomial of bi-degree $(n_1,n_2)$, such that either $n_1=-1$ or $n_2=1$. Then a Laurent polynomial $M(z)$ which is a polynomial in $z$ may not satisfy \eqref{1l+}
for $i\geq i_0,$ where $i_0\geq 0,$ unless $M(z)\equiv 0.$
\et

\pr Indeed, if $M(z)$ is a polynomial, then \eqref{lau} 
is equivalent to 
\be \la{lau1}
\sum_{i\in J_{0}}M(L^{-1}_i(t))\equiv 0.
\ee
If $n_2=1$, then \eqref{lau1} immediately implies that $M(z)\equiv 0$ since in this case $J_{0}$ contains a single element.
Suppose now that $n_1=-1$ and denote by $L^{-1}_{\infty}(z)$ a unique branch of $L^{-1}(z)$ for 
which $\lim_{z\to \infty}L^{-1}_{\infty}(z)=\infty.$
It follows from the transitivity of the monodromy group $G_L$ of $L(z)$ that
there exists $\sigma\in G_L$ such that acting on equality \eqref{lau1} by $\sigma$ 
we obtain the equality 
\be \la{lau2}
M(L^{-1}_{\infty}(t))+\sum_{i\in J_{0}\setminus j}M(L^{-1}_i(t))=0,
\ee where $j\in J_0.$
Since for any $M(z)\not \equiv 0$ we have:
$$\lim_{t\to \infty}M(L^{-1}_{\infty}(t))=\infty$$ while 
$$ \lim_{t\to \infty}M(L^{-1}_i(t))=0, \ \ \ i\in J_{0},$$ 
equality \eqref{lau2}
implies that $M(z)\equiv 0.$ 
\qed

\bt \la{d4} Let $L(z)$ be a proper Laurent polynomial of bi-degree $(n,p)$, where $p$ is a prime, and $M(z)\not \equiv 0$ be a polynomial in $z$ such that \eqref{1l+} holds 
for $i\geq i_0,$ where $i_0\geq 0.$
Then $L(z)=L_1(z^p)$ for some Laurent polynomial $L_1(z)$ while $M(z)$ is a linear 
combination of the monomials $z^{j},$ where $j \not\equiv 0 \, \mod p.$
\et

\pr 
Show first that
$J_{\infty}$ is a block of an imprimitivity system for the monodromy group $G_L$ of $L(z).$
Indeed, if $J_{\infty}$ is not a block, then there exists $\sigma\in G_L$ such that $\sigma\{J_\infty\}\cap J_\infty\neq \emptyset$
and $\sigma\{J_\infty\}\cap J_{0}\neq \emptyset$. This implies that 
acting on equality \eqref{lau1} by $\sigma$ 
we obtain the equality 
\be \la{lau3}\sum_{i\in A }M(L^{-1}_i(t))+
\sum_{i\in B}M(L^{-1}_i(t))=0,
\ee where $A$ is a subset of $J_{0}$ and $B$
is a {\it proper} subset of $J_{\infty}.$

Without loss of generality we may assume that $L^{-1}_i(t),$ $1\leq i \leq p+n,$ are numerated in such a way that $J_{\infty}=\{1,2,\dots, p\}$ and that to the loop around infinity corresponds the element 
\be \la{elem}(12\dots p)(p+1 p+2 \dots p+n) \ee
of $G_L$. Then Puiseux series of $L^{-1}_i(t),$ $1\leq i \leq p,$
at infinity are
$$ 
L^{-1}_i(z)=\sum_{k=-1}^{\infty}u_k\varepsilon_p^{(i-1)k}\left(\frac{1}{z}\right)^{\frac{k}{p}}, 
$$ where $\varepsilon_p=exp(\frac{2\pi \sqrt{-1}}{p})$ and $u_{-1}\neq 0.$ Therefore, 
$$ 
M(L^{-1}_i(z))=\beta\varepsilon_p^{(i-1)m}z^{\frac{m}{p}}+o(z^{\frac{m}{p}}), \ \ \ 
1\leq i \leq p,
$$ where $m=\deg M(z),$ $\beta=u_{-1}^{m}\neq 0.$
On the other hand, for any $i,$ $p+1\leq i\leq n+p,$ near infinity we have: 
$$M(L^{-1}_i(z))=o(1).$$ Therefore, for the coefficient $\gamma$ of $z^{\frac{m}{p}}$ in the Puiseux 
series of the function in the left side of \eqref{lau3} the equality 
\be \la{num} \gamma=\beta\sum_{j\in B}\varepsilon_p^{(j-1)m}.\ee
Thus, if we will show that $\gamma\neq 0,$ then the contradiction obtained
will imply that $J_{\infty}$ is a block.

Set
$$r(z)=\sum_{j\in B}z^{j-1}.$$ Clearly, $\gamma=\beta r(\varepsilon_p^{m}).$  
Since $p$ is a prime, the number $\varepsilon_p^{m}$ is either 1 or
a primitive $p$-th root of unity. In the first case obviously $\gamma\neq 0$.
On the other hand, in the second case the equality 
$r(\varepsilon_p^{m})=0$  
implies that the $p$-th cyclotomic polynomial $\Phi_p(z)$ divides $r(z)$ in the ring $\Z[z].$
However, this is impossible since $$\Phi_p(z)=1+z+z^2+\dots +z^{p-1},$$ while 
$B$ is a proper subset of $J_{\infty}$. Therefore, $\gamma\neq 0,$ and hence $J_{\infty}$ is a block.

Since the group $G_L$ contains a block, the Laurent polynomial $L(z)$ may be decomposed into a composition $L(z)=A(B(z))$ of rational functions $A(z)$ and $B(z)$ of degree greater
than one. 
Furthermore, since $J_{\infty}$ is a block and element \eqref{elem}
transforms $J_{\infty}$ to itself, the function  
$A(z)$ has two poles. Taking into 
account that the bi-degree of $L(z)$ is $(n,p)$, this implies that there exists a Laurent polynomial $L_1(z)$ such that $L(z)=L_1(z^p),$ where the 
bi-degree of $L_1(z)$ is $(n/p,1).$
Clearly, $M(z)$ can be written
as $$M(z)=M_1(z^p)+M_2(z),$$ where $M_1(z)$ is a polynomial in $z$ and 
$M_2(z)$ is a combination of the monomials $z^{j},$ where $j \not\equiv 0 \, \mod p.$ Furthermore, clearly $M_2(z)$ satisfies \eqref{1l+}.
Therefore, $M_1(z^p)$ also should satisfy \eqref{1l+}.
Since after the change of variable this implies that $M_1(z)$ 
satisfies \eqref{1l+} for $L(z)=L_1(z)$ it follows from  
Theorem \ref{d3} that $M_1(z)\equiv 0$. 
\qed

\section{Bautin index for the Laurent polynomial moment problem}
In this section we study the following problem: for Laurent polynomials $L(z),$ $M(z)$ to find a number $i_0$ such that the vanishing of the integrals 
\be
\int_{S^1} L^i(z)d M(z)= 0,
\ee for $i$ satisfying $0\leq i \leq i_0$ implies that they vanish for all 
$i\geq 0.$ 
A similar problem for {\it polynomials} was studied in the recent paper by V. Kisunko \cite{kis}, where a solution was given in the case 
of generic position. The approach of 
\cite{kis} is based on the fact that the function $I_{\infty}(t)=I_{\infty}(q,P,\gamma,t)$ satisfies a Fuchsian linear differential equation (see, e.g., \cite{pry}, p. 250).
In \cite{kis} this equation was calculated explicitly in the case where $P(z),q(z)$ are polynomials, and this
permitted to estimate a maximal order of a zero of $I_{\infty}(t)$ at infinity, in case where $I_{\infty}(t)\not\equiv 0$, 
via degrees of $P(z),q(z)$ implying a bound needed.

Since Theorems \ref{t1}, \ref{t2} give an explicit expression 
for the function $I_{\infty}(t)$, they provide an approach to the 
problem in the general case. 
We demonstrated this approach below in the case where 
$L(z),$ $M(z)$ are Laurent polynomials.

For a Laurent polynomial $L(z)$ of degree $n$ define numbers $f_i,$ $1\leq i \leq n,$ as follows:
$f_i,$ $1\leq i \leq n,$ equals 1 if $i\in J_0$ and 0 otherwise. Further, define a number $N(L)$ as 
the number of different vectors in the collection 
\be \la{coll} (f_{\sigma(1)},f_{\sigma(2)}, \dots f_{\sigma(n)}),\ \ 
\sigma\in  G_L.\ee Notice that obviously $N(L)\leq \vert G_L \vert\leq n!.$
Finally, for a function $\psi(t)$
whose Puiseux series 
at infinity is
\be \la{opop} \psi(t)=\sum_{k= j}^{\infty}w_k\left(\frac{1}{t}\right)^{\frac{k}{l}},\ee
where $l\geq 1$ and $w_j\neq 0$,
set $\ord_{\infty}\psi(t)=j/l.$

\bt \la{baut} Let $L(z),$ $M(z),$ $\deg L(z)=n,$ $\deg M(z)=m,$
be Laurent polynomials such that the equality 
\be\la{2l}
\int_{S^1} L^i(z)d M(z)= 0,
\ee
holds for all $i$ satisfying $0\leq i \leq m(N(L)-1)+1.$ 
Then \eqref{2l} holds for all $i\geq 0.$ In particular, equalities \eqref{2l} hold for all $i\geq 0$ whenever they hold for
all $i$ satisfying $0\leq i \leq m(n!-1)+1.$
\et 
\pr 
Set $$\psi(t)=\int I_{\infty}(t)d t=\frac{1}{2\pi i}\sum_{k=1}^{\infty}\frac{m_k}{k}\left(\frac{1}{t}\right)^{k}.$$ Clearly, we only 
must show that if \be \la{xuiak} \ord_{\infty}\psi(t)> m(N(L)-1),\ee then
$\psi(t)\equiv 0.$

Equality \eqref{eps} implies that  
\be \la{eqqq} \psi(t)= \sum_{i\in J_0} M(L^{-1}_i(t))
-\sum_{i=1}^{n}M_0(L^{-1}_i(t)).\ee 
Since $\psi(t)$ is a sum of algebraic
functions, $\psi(t)$ itself is an algebraic function and
therefore satisfies an irreducible algebraic equation
\be \la{urr}
y^N(t)+a_1(t)y^{N-1}(t)+\ ... \ + a_N(t)=0,\ \ \ a_j(t)\in \C(t),
\ee
whose roots $\psi_j(t),$ $1\leq j \leq N,$ are all possible analytic continuations of $\psi(t)$ 
and whose coefficients are elementary symmetric functions of $\psi_j(t),$ $1\leq j \leq N$. 
Furthermore, since $M(t),L(t)$ are Laurent polynomials, it follows from \eqref{eqqq} that the functions   
$\psi_j(t),$ $1\leq j \leq N,$
have no poles in $\C$ and therefore the functions $a_j(t),$ $1\leq j \leq N,$
are polynomials. Finally, 
since $$\phi(t)=\sum_{i=1}^{\deg L}M_0(L^{-1}_i(t))$$ is a rational function, the inequality 
\be \la{wewe} N\leq N(L)\ee holds.

Let $(n_1,n_2)$ be the bi-degree of $L(t).$  
The Puiseux series at infinity of branches $L^{-1}_i(t),$ $i\in J_0,$
have the form 
\be \la{aa1}
\sum_{k=1}^{\infty}v_{k,i}\left(\frac{1}{t}\right)^{\frac{k}{\ n_1}},
\ee where $v_{1,i}$ is distinct from zero,
while Puiseux series at infinity of branches $L^{-1}_i(t),$ $i\in J_{\infty},$
have the form 
\be \la{aa2}
\sum_{k=-1}^{\infty}\tilde v_{k,i}\left(\frac{1}{t}\right)^{\frac{k}{\ n_2}},
\ee where $\tilde v_{-1,i}$ is distinct from zero.
Since Puiseux series at infinity of the functions $M(L^{-1}_i(t))$, $1\leq j \leq N,$ and $M_0(L^{-1}_i(t))$, $1\leq j \leq N,$
may be obtained by the substitution of series \eqref{aa1}, \eqref{aa2} into $M(t)$ and $M_0(t)$, this implies that 
for any $j,$ $1\leq j \leq N,$ the inequality \be \la{ed} \ord_{\infty}\psi_j(t)\geq -m\ee holds.
Therefore, since $a_j(t),$ $1\leq j \leq N,$ are elementary symmetric functions of $\psi_j(t),$ $1\leq j \leq N,$ the inequalities
$$\ord_{\infty}a_j(t)\geq - mj, \ \ \ 1\leq j \leq N,$$
hold and hence 
\be \la{popo} \deg a_j(t)=-\ord_{\infty}a_j(t)\leq mj, \ \ \ 1\leq j \leq N.\ee

Now we are ready to show that if \eqref{xuiak} holds, then 
$\psi(t)\equiv 0.$ Indeed, assume the inverse. Then for the coefficient $a_N(t)$ in \eqref{urr} the inequality $\ord_{\infty}a_N(t)\leq 0$ holds. On the other hand, \eqref{xuiak} and \eqref{wewe} imply that $\ord_{\infty}\psi(t)> m(N-1)$
and hence for any $i,$ $1\leq i\leq N,$ taking into account inequalities \eqref{popo}, 
we have:
\begin{multline*} \ord_{\infty}\{a_{N-i}(t)\psi^i(t)\}\geq \ord_{\infty}\{a_{N-i}(t)\psi(t)\}= \\
=\ord_{\infty} \psi(t)-\deg a_{N-i}(t)\geq \ord_{\infty} \psi(t)-m(N-1)>0\end{multline*}
(here we set $a_0(t)\equiv 1$).  
Therefore, 
$$\ord_{\infty} \{\psi^N(t)+a_1(t)\psi^{N-1}(t)+\ ... \ + a_{N-1}(t)\psi(t) \}>0$$  
in contradiction with 
$$\psi^N(t)+a_1(t)\psi^{N-1}(t)+\ ... \ + a_N(t)=0\ \ \ \ \Box$$

\noindent{\bf Remark.} The proof of Theorem \ref{baut}
uses the same ideas as Section 2.4 of \cite{pp}.
Notice that corresponding formulas in \cite{pp} on the page 758 contain misprints. Namely, all printed powers of the expression $m/n$ are actually its 
factors.

\section{\la{cri} Relations between the rationality and the reducibility of $I_{\infty}(q,P,\gamma,t)$} 
\subsection{\la{inf} Definition of the subspace $M_{P,\gamma}$} 
Let $\gamma$ be a curve and $P(z)$ be a rational function such that poles of $P(z)$ do not lie on one side of $\gamma$.
By Corollary \ref{gen+}, if $P(z)$ is in generic position, then the rationality of $I_{\infty}(q,P,\gamma,t)$ for a rational function $q(z)$ implies that $q(z)=\tilde q(P(z))$ for some rational function 
$\tilde q(z).$ For arbitrary $P(z)$ such a statement fails to be true.
However, in many cases the rationality of $I_{\infty}(q,P,\gamma,t)$
still implies that $I_{\infty}(q,P,\gamma,t)$ is reducible. 

For example, if $\gamma$ is a non-closed curve such that its end points $a,b$ are not ramification points of $P(z)$, then the rationality of 
$I_{\infty}(q,P,\gamma,t)$ implies its reducibility. Indeed, by Corollary 
\ref{rysa} in this case $P(a)=P(b)$, and equality \eqref{fo3} from Proposition \ref{egg} reduces 
to the equality \be \la{eby} \left(\frac{q}{P'}\right) (P^{-1}_{i_1}(z))=\left(\frac{q}{P'}\right)(P^{-1}_{i_2}(z))\ee
for some $i_1\neq i_2,$ $1\leq i_1,i_2\leq n.$ 
Therefore, by Lemma \ref{compos}  
there exist rational functions $R(z),$ $\tilde P(z)$, and 
$W(z)$ with $\deg W(z)>1$ such that
\be \la{iilo} \left(\frac{q}{P'}\right) (z)=R(W(z)), \ \ \ P(z)=\tilde P(W(z))\ee and hence 
\eqref{poi} holds for $\tilde q(z)=R(z)\tilde P'(z)$ since \eqref{iilo} yields that 
$$q(z)=R(W(z))P^{\prime}(z)=R(W(z))\tilde P'(W(z))W'(z)=\tilde q(W(z))W'(z).$$

In this section we in a sense describe the class of pairs $P(z),$ $\gamma$ for which the rationality of $I_{\infty}(q,P,\gamma,t)$ implies its reducibility.
For given $P(z)$ and $\gamma$, such that poles of $P(z)$ do not lie on one side of $\gamma$,
a natural necessary condition for the existence of $q(z)$ such that
$I_{\infty}(q,P,\gamma,t)$ is rational but is not reducible
may be formulated as follows. Let 
$$\sum_{i=1}^nf_{s,i}\left(\frac{q}{P'}\right)(P^{-1}_{i}(z))=0,\ \ \ 1\leq s \leq k,$$ 
be the system of equations from Theorem \ref{t3} and let $M_{P,\gamma}$ be 
a linear subspace of $\mathbb Q^n$ 
generated by the vectors 
$$(f_{s,\sigma(1)}, f_{s,\sigma(2)},\, ...\,, f_{s,\sigma(n)}), \ \ \ \sigma\in G_P, \ \ \ 1\leq s \leq k,$$
where $G_P$ is the monodromy group of $P(z)$ and $n=\deg P(z).$ 
By Corollary \ref{polusa} the subspace $M_{P,\gamma}$ is not zero-dimensional.
Furthermore, by  construction, $M_{P,\gamma}$ is invariant with respect to the permutation representation of $G_P$ over $\Q$ and it follows from Theorem \ref{t3} by the analytic continuation, that for any vector $\vec v\in M_{P,\gamma}$, 
$\vec v=(v_1,v_2,\dots, v_n)$, the equality 
\be \la{vvbbcc} \sum_{i=1}^nv_{i}\left(\frac{q}{P'}\right)(P^{-1}_{i}(z))=0\ee holds.

Observe now that if $M_{P,\gamma}$ contains
a vector of the form 
$\vec{e}_i-\vec{e}_j,$  $i\neq j,$
where $\vec e_i,$ $1\leq i \leq n,$ denote vectors of the Euclidean basis of $\Q^n,$ 
then the rationality of $I_{\infty}(q,P,\gamma,t)$ implies its
reducibility since for such a vector 
equality \eqref{vvbbcc} implies \eqref{eby} and \eqref{iilo}.
Therefore, a necessary condition for the existence of 
$q(z)$ such that $I_{\infty}(q,P,\gamma,t)$ is rational but is not reducible
is that $M_{P,\gamma}$
contains no vectors of the form $\vec{e}_i-\vec{e}_j,$ $i\neq j,$ and in this section, using a general result of 
\cite{gi2}, we prove (Theorem \ref{gir2} below)
that this condition is also sufficient.
As an application 
we show that the requirement of  
Theorem \ref{gen} can not be weakened to the requirement of  indecomposability of $P(z)$ already for Laurent polynomials.

\subsection{Girstmair's theorem}
Let $f(t)\in K[t]$ be an irreducible polynomial over a field of characteristic zero $K$. Denote by $x_1,x_2, \dots, x_n$ roots of $f(t)$, by $L$ the field $K(x_1,x_2, \dots, x_n)$, and by $G$ the Galois group $\Gal(L/K).$ Usually we will identify $G$ with a permutation group acting on 
the set $\{1,2,\dots,n\}$ setting 
$\sigma(i)=j$ 
if $x_j=\sigma(x_i),$ $1\leq i,j \leq n$.
In this subsection, following \cite{gi2}, we sketch a solution of the following problem: {\it under what conditions on a collection $W$ of vectors from $K^n$
there exists a rational function 
$R(t)\in K(t)$ 
such that 
for all $\vec w\in W,$ $\vec w=(w_1,w_2, \dots, w_n),$ the equality 
\be \la{gg1} w_1R(x_1)+w_2R(x_2)+\dots +w_nR(x_n)=0\ee holds, and 
$R(x_i)\neq R(x_j)$ for any $i\neq j,$ $1\leq i,j\leq n.$}
If such a function $R(t)$ exists we will say that $W$ is {\it admissible}.

Notice that if \eqref{gg1} holds for vectors $\vec w_1, \vec w_2$, then it holds for the vector $a\vec w_1+b\vec w_2$, $a,b\in K.$  
Furthermore, for any element $\sigma\in G$, acting on equality \eqref{gg1} by
$\sigma$ and replacing $\sigma$ by $\sigma^{-1}$, we obtain the equality 
$$ w_{\sigma(1)}R(x_1)+w_{\sigma(2)}R(x_2)+\dots +w_{\sigma(n)}R(x_n)=0.$$ 
Therefore, equality \eqref{gg1} holds for all $\vec w\in W$ if and only if it holds for all vectors from the linear subspace 
of $K^n$ generated by the vectors 
$$\vec w^{\sigma}=(w_{\sigma(1)},w_{\sigma(2)}, \dots w_{\sigma(n)}),\ \ 
w\in W, \ \ 
\sigma\in  G.$$ Thus, without loss of generality we may assume 
that the collection $W$ is a linear subspace of $K^n$ invariant with respect to 
the permutation representation of $G$ on $K^n.$

We start from reformulating the problem above in the form it was considered in \cite{gi2}. Fix a root $x$ of $f(t)$. Denote by $H$ the stabilizer $G_{x}$ of $x$ in $G$ and by $G/H=\{\bar s\,:\, s\in G\}$ the set of left cosets $\bar s=sH$ of the subgroup $H$ in $G$. Further, denote by $K[G]$ the group ring of $G$ over $K$ and by $K[G/H]$ a $K$-module with the basis $(\bar s\,:\, \bar s\in G/H)$. 
Thus, elements of $K[G]$ have the form 
\be \la{form} \lambda= \sum_{s\in G}l_ss,\ \ \ l_s\in K,\ee
while elements of 
$K[G/H]$ have the form
\be \la{al} \alpha=\sum_{\bar s\in G/H} a_{\bar s}\bar s, \ \ a_{\bar s}\in K.\ee  
Notice that $K[G/H]$ is a $K[G]$-module with respect to the scalar multiplication defined by the formula
$$g\bar s= \overline{gs}, \ \ \ g\in G, \ \ \ \bar s \in G/H.$$

If $y\in L$ satisfies $G_y=H$, then for and any $\alpha\in K[G/H]$ 
defined by \eqref{al} the expression 
$$\alpha y= \sum_{\bar s\in G/H}  a_{\bar s}sx$$ is a well defined element of $L.$ We say that a subset $M$ of $K[G/H]$ is {\it admissible} if there exists $y\in L$ such that $G_y=H$ and for 
any $\alpha\in M$ the equality $\alpha y=0$ holds. Clearly, 
$M$ is admissible if and only if the $K[G]$-submodule of $K[G/H]$
generated by $M$ is 
admissible so without loss of generality we may assume that $M$ is 
such a submodule.

Recall that linear subspaces of $K^n$, invariant with respect to 
the permutation representation of $G$ on $K^n$, are in one-to-one correspondence with $K[G]$-submodules of $K[G/H].$ Namely, to a subspace 
$W$ corresponds a submodule $\widehat W$ consisting of elements 
$$\alpha=\sum_{i=1}^n a_{i}\bar s_i,$$
where $s_i,$ $1\leq i \leq n,$ is an element of $G$ which transforms
$1$ to $i$ and 
$\vec a=(a_1,a_2,\dots, a_n)$ runs elements of $W$.

\bp A linear subspace $W$ of $K^n$, invariant with respect to 
the permutation representation of $G$ on $K^n$, is admissible if and only 
if the corresponding $K[G]$-submodule $\widehat W$ of $K[G/H]$ is admissible.
\ep
\pr Show first that if $W$ is admissible, then we can set $y=R(x_1).$ Indeed, since 
$y\in K(x_1)$ we have $H\subseteq G_y$. Furthermore, $H$ may not be a proper subgroup of $G_y$ since otherwise the length of the orbit of 
$y$ under the action of $G$ would be strictly less than $n$ in contradiction with the conditions that all $R(x_i)$, $1\leq i\leq n,$ 
are different between themselves.
Finally,
$\alpha y=0$ for any $\alpha \in \widehat W.$ 

In other direction, if $\widehat W$ is admissible and $y$ is an element of $L$ such that $\alpha y=0$ for all $\alpha\in \widehat W$, then $G_y=H$ implies that $y\in K(x_1).$ 
Therefore, there exists $R(t)\in K(t)$ such that $y=R(x_1)$ and 
for such $R(z)$ equality \eqref{gg1} holds  
for all $\vec w\in W.$
Furthermore, since the length of the orbit of 
$y$ under the action of $G$ equals $n$, all $R(x_i)$, $1\leq i\leq n,$ 
are different between themselves.  \qed

\bt[\cite{gi2}] \la{gir} A $K[G]$-submodule $\widehat W$ of $K[G/H]$ is admissible if and only if $\widehat W$ contains no elements   
$\bar s_1- \bar s_2,$ $s_1,s_2\in G,$ unless $\bar s_1= \bar s_2.$ 
Equivalently, 
a linear subspace $W$ of $K^n$, invariant with respect to 
the permutation representation of $G$ on $K^n$, is admissible if and only 
$W$ contains no vectors $\vec e_i-\vec e_j$, $1\leq i,j\leq n,$
unless $i= j.$ 
\et
\pr Indeed, if $W$ contains a vector $\vec w=\vec e_i-\vec e_j$, $i\neq j,$ $1\leq i,j\leq n,$ then equality \eqref{gg1} implies that $R(x_i)=R(x_j).$

In other direction, assume that $\widehat W$ contains no elements $\bar s_1- \bar s_2,$ $s_1,s_2\in G,$ such that $\bar s_1\neq \bar s_2.$ Consider the canonical $K[G]$-linear map 
$$\rho\,:\, K[G]\rightarrow K[G/H]$$ which maps $s$ to $\bar s$, 
and let $\gamma=\rho^{-1}(\widehat W)$
be the inverse image of $\widehat W.$
Since $K[G]$ is semisimple, the ideal $\gamma$ is generated by an idempotent element $\v$. Notice that for any $\lambda\in \gamma$ the equalities $\lambda=a\v,$ $a\in K[G],$ and $\v^2=\v$ imply that   
$\lambda \v=\lambda.$
Set $\mu=1-\v$. Then for any $\alpha\in \widehat W$ the equality $\alpha \mu=0$ holds. Indeed, if $\lambda$ is an element of 
$\gamma$ such that $\rho(\lambda)=\alpha$ then we have 
$$\alpha\mu=\lambda\mu=\lambda-\lambda\v=0.$$

Furthermore, for any $s\in H$ the element $s - 1$ is in the kernel of $\rho$ and therefore in $\gamma$. Hence $s-1=(s-1)\v$ and
$(s-1)\mu=0$ implying $H\subseteq G_{\mu}.$ On the other hand, for any $s\in G_{\mu}$ we have $(s-1)\mu=0.$ Therefore, $s-1=(s-1)\v$ and 
$s-1\in \gamma.$ This implies that $\bar s -\bar 1$ is in $\widehat W$
and therefore $s\in H$ by the assumption. This proves that $H= G_{\mu}.$

Finally, let us show that from the existence $\mu \in  K[G]$, such that $H= G_{\mu}$ and for any $\alpha\in \widehat W$ the equality $\alpha \mu=0$ holds, follows that $\widehat W$ is admissible. For this purpose observe that by the normal basis theorem there exists an 
element $x\in L$ such that $gx,$ $g\in G$, is a basis of $L$ over $K$.
Set now $y=\mu x$. Then obviously for any $\alpha\in \widehat W$ the equality $\alpha y=0$ holds and $H\subseteq G_y$. Furthermore, $H= G_y$. Indeed,
if there exists $g_0\in G$ such that $g_0y=y$ but $g_0\mu \neq \mu$, then it follows from equalities 
$g_0\mu x=y,$  $\mu x=y$ that $gx,$ $g\in G,$ are linearly dependent over $K.$
\qed

\subsection{Existence of $q(z)$ with rational but not reducible $I_{\infty}(q,P,\gamma,t)$} 
Theorem \ref{gir} permits to solve the problem posed in Subsection \ref{inf} as follows.

\bt\la{gir2} Let $\gamma$ be a curve and $P(z)$ be a rational function of degree $n$ such that poles of $P(z)$ do not lie on one side of $\gamma$.
Then a rational function $q(z)$, such that 
$I_{\infty}(q,P,\gamma,t)$ is rational but is not reducible, exists if and 
only if the subspace $M_{P,\gamma}$ contains no vectors $\vec e_i-\vec e_j$, $1\leq i,j\leq n,$ unless $i= j.$ 
\et
\pr As it was already observed in Subsection \ref{inf} the requirement of the theorem is necessary. On the other hand, 
since vectors with rational coefficients which are linear independent over $\Q$ remain linearly independent over $\C(z)$, if this requirement is satisfied, then the subspace $\widehat M_{P,\gamma}$
of $(\mathbb C(z))^n$, 
generated over $\C(z)$ by the same vectors 
$$(f_{s,\sigma(1)}, f_{s,\sigma(2)},\, ...\,, f_{s,\sigma(n)}), \ \ \ \sigma\in G_P, \ \ \ 1\leq s \leq k,$$
which generate $M_{P,\gamma}$ over $\Q$, still contains no vectors $\vec e_i-\vec e_j$, $1\leq i,j\leq n,$ unless $i= j.$ 

Therefore, applying
Theorem \ref{gir} to the roots 
$$x_1=P^{-1}_1(z),\ \ x_2=P^{-1}_2(z),\ \  \dots, \ \ x_n=P^{-1}_n(z)$$ of the 
polynomial $P(x)-z=0$ over the field $\C(z)$, we conclude that there exists a rational function $\widehat R(t) \in \C(z)(t)$ (whose coefficients are rational functions !) such that
for any vector $\vec v\in \widehat M_{P,\gamma}$, 
$\vec v=(v_1,v_2,\dots, v_n)$, 
the equality 
\be \la{vvbbcc2} \sum_{i=1}^nv_{i}\widehat R(P^{-1}_{i}(z))=0\ee holds
and all $\widehat R(P^{-1}_i(z))$, $1\leq i\leq n,$ 
are different between themselves. Furthermore, since we can write $z$ as $z=P(P^{-1}_{i}(z)),$
$1\leq i\leq n,$ there exists a polynomial $R(t)\in \C(t)$ 
(whose coefficients now are just complex numbers) such that 
$$R(P^{-1}_i(z))=\widehat R(P^{-1}_i(z)), \ \  \ 1\leq i\leq n.$$
Setting now $q(z)=R(z)P'(z)$ we see that for any vector $\vec v\in M_{P,\gamma}$ the equality \eqref{vvbbcc} holds. Furthermore, equality \eqref{poi} is impossible since otherwise Lemma \ref{compos} would imply 
that $$R(P^{-1}_{i}(z))=R(P^{-1}_{j}(z))$$ for some $i\neq j,$ $1\leq i,j \leq n.$ \qed

Using Theorem \ref{gir2} one can prove the existence of an indecomposable Laurent polynomial $L(z)$ for which there exists a rational function $q(z)$, such that $I_{\infty}(q,P,\gamma,t)$ is rational but is not reducible, without an actual calculation 
$L(z)$ and $q(z).$ Indeed,  
let $L(z)$ be a Laurent polynomial whose constellation 
is shown on Fig. 4 and 
\begin{figure}[ht]
\epsfxsize=5.5truecm
\centerline{\epsffile{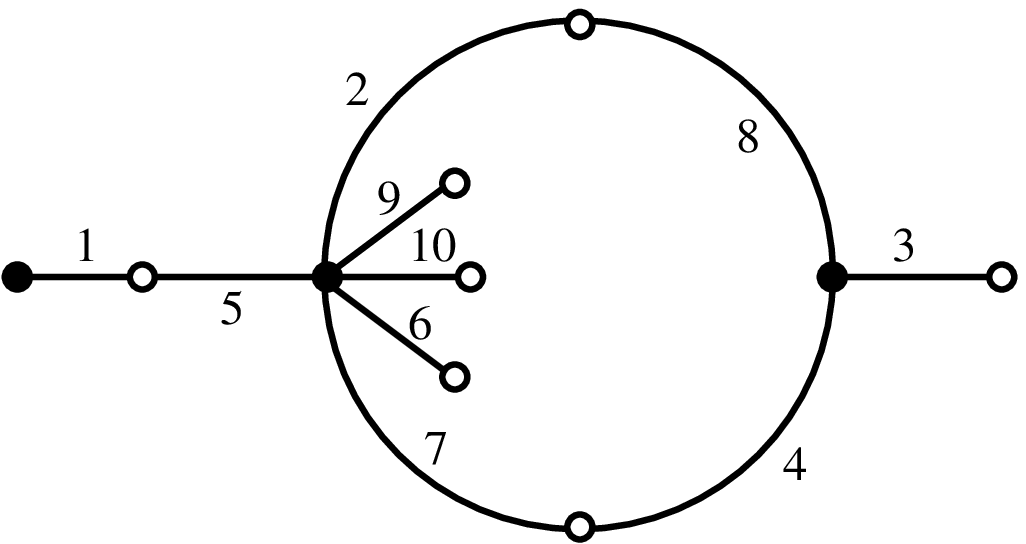}}
\smallskip
\centerline{Figure 4}
\end{figure}
whose monodromy group $G$ is generated by the permutations
$\alpha=(2,5,7,6,10,9)(3,8,4)$ and $\beta=(1,5)(2,8)(4,7)$ 
(since $L(z)$ has only two finite critical values, 
in correspondence with the notation of ``dessins d'enfants'' theory, we picture here $1$-vertices 
as ``black'', $2$-vertices 
as ``white'', and do not mark non-numerated vertices at all).
This choice of $L(z)$ is motivated by the fact that the action of $G$ on branches of $L^{-1}(z)$ is permutation equivalent to the action of the group $S_5$ on two element subsets of $\{1,2,3,4,5\}$. Since it is well known that the last action is primitive while the corresponding matrix representation of dimension $10$ over $\Q$ is not a sum of $V_{\Q}$ and 
$V_{\Q}^{\perp},$ one can expect that $L(z)$ provides
a desired example.

By Theorem \ref{t3} the function 
$I_{\infty}(q,P,\gamma,t)$ is rational if and only if the equality 
\be\la{rys}
Q(L_{2}^{-1}(z)) \,-\, Q(L_{7}^{-1}(z)) \,+\,
Q(L_{4}^{-1}(z)) \,-\, Q(L_{8}^{-1}(z)) \,\equiv\, 0\,,
\ee
holds. Therefore, the subspace $M_{L,S^1}$ is generated by the 
single vector $$\vec v=(0,1,0,1,0,0,-1,-1,0,0).$$ Show now that $M_{L,S^1}$ may not contain a vector $w$ of the form \be \la{vvee} w=\vec e_i-\vec e_j, \  \ i\neq j, \ \ 1\leq i,j\leq n.\ee 

Consider the vector subspace $V$ of $\Q^{10}$ generated by the vectors
\begin{eqnarray*}
\vec v_1 & = & (1,0,0,0,1,1,0,0,1,0)\,, \\
\vec v_2 & = & (1,1,0,0,0,0,1,0,0,1)\,, \\
\vec v_3 & = & (0,1,1,0,0,1,0,1,0,0)\,, \\
\vec v_4 & = & (0,0,1,1,0,0,1,0,1,0)\,, \\
\vec v_5 & = & (0,0,0,1,1,0,0,1,0,1).
\end{eqnarray*}
Since $\alpha$ and $\beta$ permute the vectors $\vec v_i,$ $1\leq i \leq 5,$ between themselves, $V$ is $\rho_G$-invariant. 
Furthermore, 
since $\vec v$ is orthogonal to $\vec v_i,$ $1\leq i \leq 5,$
the inclusion $$M_{L,S^1}\subseteq V^{\perp}$$ holds. On the other hand, it is easy to see that for any vector \eqref{vvee} there exists $\vec v_i,$ $1\leq i \leq 5,$ such that $(w,v_i)\neq 0$. Indeed, since $G$ is transitive and permute $\vec v_i,$ $1\leq i \leq 5,$ it is enough to verify this property only for $w$ whose first coordinate equals 1 and for such $w$ we may take one of the vectors $v_1$, $v_2$. 
Therefore, $M_{L,S^1}$ may not contain $w$ and hence by Theorem \ref{gir2} there exists a rational function $q(z)$ such that $I_{\infty}(q,P,\gamma,t)$ is rational but is not reducible.

For a comprehensive study of the above example in the context of the Laurent 
polynomial moment problem we refer the reader to the paper \cite{ppz}.

\vskip 0.2cm
\noindent{\bf Acknowledgments}.  The results of this paper were obtained mostly during the visit of the author to the Max-Planck-Institut f\"ur Mathematik in Spring 2009 and the author would 
like to thank the Max-Planck-Institut for the support and the hospitality.

\end{document}